\documentclass[12pt]{amsart}
\usepackage{amssymb,latexsym}
\usepackage[mathscr]{eucal}
\usepackage{graphics}
\usepackage{ulem} 
\usepackage{verbatim}

\theoremstyle{plain}
\newtheorem{theorem}{Theorem}[section]
\newtheorem{corollary}[theorem]{Corollary}
\newtheorem{lemma}[theorem]{Lemma}
\newtheorem{proposition}[theorem]{Proposition}
\newtheorem{remark}[theorem]{Remark}
\newtheorem{definition}[theorem]{Definition}
\newtheorem{example}[theorem]{Example}
\newtheorem{question}[theorem]{Question}
\def\sideremark#1{\ifvmode\leavevmode\fi\vadjust{\vbox to0pt{\vss
\hbox to 0pt{\hskip\hsize\hskip1em
\vbox{\hsize2cm\tiny\raggedright\pretolerance10000
\noindent#1\hfill}\hss}\vbox to8pt{\vfil}\vss}}}

\newcommand{\be}{\begin{equation}\label}
\newcommand{\ee}{\end{equation}}
\newcommand{\bq}{\begin{equation*}}
\newcommand{\eq}{\end{equation*}}
\newcommand{\ba}{\begin{align*}}
\newcommand{\ea}{\end{align*}}
\newcommand{\bp}{\begin{proof}}
\newcommand{\ep}{\end{proof}}
\newcommand{\bL}{\begin{lemma}\label}
\newcommand{\eL}{\end{lemma}}
\newcommand{\bP}{\begin{proposition}\label}
\newcommand{\eP}{\end{proposition}}
\newcommand{\bC}{\begin{corollary}\label}
\newcommand{\eC}{\end{corollary}}
\newcommand{\bT}{\begin{theorem}\label}
\newcommand{\eT}{\end{theorem}}
\newcommand{\bR}{\begin{remark}\label}
\newcommand{\eR}{\end{remark}}
\newcommand{\bD}{\begin{definition}\label}
\newcommand{\eD}{\end{definition}}

\newcommand{\tr}{\text{Tr}~}

\newcommand{\2}{\frac{1}{2}}
\newcommand{\card}{\text{card}}

\begin{document}

\title{Strong sums of projections   in von Neumann  factors.}
\author{Victor Kaftal}
\address{University of Cincinnati\\
          Department of Mathematics\\
          Cincinnati, OH, 45221-0025\\
          USA}
\email{victor.kaftal@math.uc.edu}
\author{Ping Wong Ng}
\address{University of Louisiana\\
          Department of Mathematics\\
          Lafayette, LA, 70504\\
          USA}
\email{png@louisiana.edu}
\author{Shuang Zhang}
\address{University of Cincinnati\\
          Department of Mathematics\\
          Cincinnati, OH, 45221-0025\\
          USA}
   \email{zhangs@math.uc.edu}
\date{\today}

\begin{abstract}
This paper presents necessary and sufficient conditions for a positive bounded operator on
 a separable Hilbert space  to be the sum of a finite or infinite collection of projections
 (not necessarily mutually orthogonal), with the sum converging in the strong
 operator topology if the collection is infinite.  A similar necessary condition
  is given  when the operator and the projections are taken in
   a type II von Neumann factor, and the condition is proven to be also sufficient
    if the operator is ``diagonalizable". A simpler necessary and sufficient condition
    is given in the type III factor case.
\end{abstract}
\maketitle
\section{Introduction}

Which positive bounded operators on a separable Hilbert space can be written
as sums of projections? For finite sums, Fillmore asked this
question  and obtained the characterizations of the finite rank
operators that are sums of projections  \cite [Theorems 1] {Fp69}
(see  Corollary \ref {C:2.5}   below) and of the bounded
operators that are sums of two projections \cite [Theorems 2] {Fp69}
(see Proposition \ref {P:2.10}  below.)

For infinite sums with convergence in the strong operator topology,  this question arose naturally from work on ellipsoidal tight frames by  Dykema,  Freeman, Kornelson, Larson,
Ordower, and Weber in \cite {DFKLOW}. They proved that a sufficient
condition for a positive bounded operator $A\in B(H)^+$ to be the
sum of projections  is that its essential norm $||A||_{e} $ is larger than one (\cite [Theorem 2]
{DFKLOW}). This result served as a basis for further work by
Kornelson and Larson \cite {KkLd04} and then by Antezana,   Massey,
Ruiz, and Stojanoff \cite {AMRS} on the decomposition of positive operators
into (strongly
 converging) sums of rank-one positive operators of preset norms.

The same question can be asked relative to a von Neumann algebra $M$.
 We say that an operator $A\in M^+$  is a strong sum of projections if
 there  exists a collection of (not necessarily mutually orthogonal
 or commuting) projections $P_j\in M$ with cardinality $N\le \infty$,
 for which $A=\sum_{j=1}^N P_j$ and the series converges in the strong
  operator topology if  $N=\infty$.   The main goal of this article is to answer the
   question of which operators are strong sums of  projections.

To simplify the treatment, we consider only von Neumann factors, and we
 further assume that they are $\sigma$-finite (i.e., countably decomposable)
 so that all infinite projections are equivalent.
Thus let $H$ be a complex infinite dimensional Hilbert space and
$M\subset B(H)$ be a  $\sigma$-finite von Neumann factor. If $M$ is  
of type I, we will identify it with $B(H)$ (hence we will assume that $H$
is separable), and denote by $\tr$ the usual  normalized trace
such that $\tr P=1$ for any rank-one projection $P$. If $M$ is of
type II, $\tau$ will denote the faithful positive semifinite normal
trace, unique up to scalar multiples in the type II$_\infty$ case
and normalized by  $\tau(I)=1$ in the type II$_1$ case.  If $M$ is
only assumed to be semifinite, i.e.,  it is of type I or type II
unless specified, we will generically denote its trace by $\tau$.

The conditions  for $A$ to be a strong sum of projections are
 expressed in terms of the \textit{excess} and the \textit{defect}  parts
 of  $A$. Given  $A\in M^+$, we denote by

\ba
\quad\quad\quad\quad\quad  &\chi_A &&\text{the spectral measure of $A$}
\quad\quad\quad\quad\quad\quad \quad\quad\\
&R_A=\chi_A(0, ||A||] &&\text{the range projection of $A$}\\
&A_+:=(A-I)\chi_A(1, ||A||] && \text{the excess part}\text{ of $A$} \\
&A_-:= (I-A)\chi_A(0,1) && \text{the defect part} \text{ of $A$}.
\end{align*}
Thus we have the decomposition

\be{e:1}
A= A_+ - A_-+R_A.
\ee

A positive operator $A$ is said to be {\sl diagonalizable} if
$A=\sum \gamma_j E_j$ for some  $\gamma_j> 0$
 and  mutually orthogonal projections $\{E_j\}$  in $M$. Diagonalizable
 operators are  also called \textit{discrete} and are the most accessible operators in a type
  II  factor (e.g., see  \cite {AmMpCarp}.)

The main results of this article are collected in the following
theorem.

\bT{T:1.1} Assume that $M$ is a $\sigma$-finite von Neumann
factor and $A\in M^+$.
\item[(i)] Let  $M$ be of type I. Then $A$ is a strong sum of projections
 if and only if either $\tr(A_+)=\infty$ or $\tr(A_-) \le \tr(A_+) < \infty$
 and $\tr(A_+) - \tr(A_-)\in \mathbb N\cup \{0\}$.(Theorems \ref {T:6.6},
 \ref {T:4.3}, and \ref {T:3.3}.)
\item [(ii)]  Let $M$ be of type II and $A$ be diagonalizable. Then $A$
is a strong sum of projections if and only if $\tau(A_+)\ge
\tau(A_-)$. The condition is necessary even when $A$ is not
diagonalizable. (Theorems \ref {T:6.6}, \ref {T:5.2}, and
\ref {T:3.3}.)
\item[(iii)] Let $M$  be of type III. Then  $A$ is a strong sum of projections
if and only if either $||A||>1$ or $A$ is a projection. (Corollary
\ref {C:6.4} and Theorem \ref {T:3.3}.) \eT

 \bR {R:1.2}
The statement (i) above extends the  sufficient condition obtained
in
 \cite [Theorem 2] {DFKLOW}. In fact,  it is elementary to show that
$||A||_{e}>1 $ implies that $ \tr(A_+)=\infty; $ however, the reverse
implication is false.

\eR

The necessary conditions in Theorem \ref {T:1.1} are obtained
via the
 frame theory type construction  of Proposition \ref {P3.1}  that links decomposability of an
operator $A$ into a strong sum of projections to the condition that the  identity is 
the  ``diagonal" of $W^*AW$ for some  partial  isometry $W$ with
$W^*W=R_A$. For instance, the integrality condition in the $B(H)$
case (Theorem \ref {T:1.1} (i)) when $\tr(A_+)< \infty$ emerges
naturally from the fact that $\tr(A_+) - \tr(A_-)$ coincides with
the trace of the projection  $I-WW^*$.

A modification of these arguments  provides an alternative proof of the necessity of the ``integrality condition" for diagonals of projections  in Kadison's   \cite [Theorem 15]{Kr02b} that identifies explicitly the integer as the difference of traces of two projections (Corollary \ref {C:3.6}.)

The basic tool for all the sufficient conditions is provided by a
$2\times2$ matrix construction that decomposes certain
diagonal matrices into the sum of a projection and a rank-one
``remainder" (Lemma \ref {L:2.1} ). This lemma serves also 
several other purposes: when applied to finite matrices it provides
in Corollary \ref {C:2.5}  another proof of Fillmore's characterization of finite sums of projections \cite [Theorem~1]{Fp69}. It can be applied  to (finite) sums of scalar multiples of
mutually orthogonal equivalent projections in a $C^*$-algebra (Lemma
\ref {L:2.6} ). It  also provides in the von Neumann algebra
setting a short constructive proof (Proposition \ref {P:2.10} )  of Fillmore's   characterization of sums of two projections \cite [Theorem 2, Corollary] {Fp69}. 

As the results of \cite {DFKLOW} suggest, the most tractable case
 is the ``infinite" one. The key special case (Lemma \ref {L:6.1})
 is when $A$  is an infinite sum of scalar multiples of mutually orthogonal
 equivalent projections in $M$ and the sum of the coefficients in the
 corresponding expansion of $A_+$ diverges.  Based on this lemma we obtain
  the sufficiency in Theorem \ref {T:1.1} for  part (iii), for part (i)
   when  $\tr(A_+)=\infty$,  and for part (ii) when  $\tau(A_+)=\infty$.

For the more delicate  ``finite trace" case in $B(H)$,  i.e., when $\tr(A_+)<\infty$,
we diagonalize $A_+$ and $A_-$ and then apply iteratively Lemma \ref {L:2.1} ,
which provides canonically a sequence of projections. The  strong convergence
 of the series of these projections is proven by reducing the problem to a
 finite dimensional construction and to three infinite dimensional special cases
 (Lemmas \ref {L:2.3}  , \ref {L:4.1}, \ref {L:4.2}, and
  Theorem \ref  {T:4.3}.)

When $M$ is of type II, Lemma \ref {L:2.1} can also be applied to diagonalizable
 operators, where the strong convergence of the ``remainders"  is obtained by
  showing that they converge in the trace-norm (Lemma \ref {L:5.1}).  Example
  \ref  {E:5.3} exhibits a non-diagonalizable operator  that is the
   sum of two projections. It  remains open whether the condition
   $\tau(A_+)\ge \tau(A_-)$ is always sufficient for $A$ to  be the strong
   sum of projections.

Von Neumann algebras are by no means the only setting in  which
 positive operators may be decomposed into sums of projections. In a separate
 paper  (\cite{KNZ08}), we will investigate the same problem for positive operators in the multiplier
 algebra $M(\mathscr A\otimes K)$ where $\mathscr A$ is a $\sigma-$unital
purely infinite simple $C^*$-algebra.

 The first and second named authors were participants in the  NSF supported
 Workshop in Linear Analysis and Probability, Texas A\&M University, 2006,
 where they first heard from David Larson about the results in \cite {DFKLOW}
 and  \cite {KkLd04} that stimulated this project.
 
\section{The matrix construction}

We start with a simple lemma which will be used in   our key
constructions.

\bL{L:2.1} Let $e$ and $f$ be two orthogonal unit vectors in $H$. For every
$\mu \ge 0$ and $ 0\le\lambda\le 1$, let
\be{e:2}
\nu: =\begin{cases} \frac{(1-\lambda)\lambda}{(1+\mu-\lambda)(\mu+\lambda)}
\quad &\text{for }\mu\ne 0\\
1 & \text{for }\mu= 0
\end{cases}\quad\quad \text{and}\quad\quad
 \rho:=\begin{cases} \frac{(1-\lambda)\mu}{\mu+\lambda}\quad &\text{for }\mu\ne 0\\
 0& \text{for }\mu= 0
\end{cases}
\ee
and let
\be{e:3}
w:=\sqrt{\rho}f- \sqrt{1-\rho}\,e \quad {and}\quad   v:= \sqrt{\nu}f+ \sqrt{1-\nu}\,e.
\ee
 Then $w\otimes w$ and $v\otimes v$ are rank-one projections and
\be{e:4}
(1+\mu)(e\otimes e) + (1-\lambda)(f\otimes f) =w\otimes w+ (1+\mu-\lambda)(v\otimes v).
\ee
\eL
\bp
It is immediate to verify that $0 \le \nu, \rho \le1$,  $w$ and $v$ are
unit vectors, and hence $w\otimes w$, $v\otimes v$ are rank-one projections
with range contained in span$\{e,\,f\}. $ Their matrix representations with respect to the basis $\{e, f\}$ are, respectively,
\[
\begin{pmatrix}
1-\rho&-\sqrt{\rho(1-\rho)}\\
-\sqrt{\rho(1-\rho)} & \rho
\end{pmatrix}
\quad \text {and} \quad
\begin{pmatrix}
1-\nu&\sqrt{\nu(1-\nu)}\\
\sqrt{\nu(1-\nu)} & \nu
\end{pmatrix}.
\]

An elementary computation shows that

\[
\begin{pmatrix}
1+\mu&0\\
0&1-\lambda
\end{pmatrix} \\
=\begin{pmatrix}
1-\rho&-\sqrt{\rho(1-\rho)}\\
-\sqrt{\rho(1-\rho)} & \rho
\end{pmatrix}
+(1+\mu-\lambda)
\begin{pmatrix}
1-\nu&\sqrt{\nu(1-\nu)}\\
\sqrt{\nu(1-\nu)} & \nu
\end{pmatrix}
\]
and hence (\ref {e:4}) holds.
\ep

If we do not require the orthogonality of the vectors $e$ and $f$, we still obtain the decomposition  in \eqref {e:4}, but the vectors $w$ and $v$ are no longer obtained as simply as in \eqref {e:3}.
With a slight generalization and a reformulation in terms of rank-one projections, we have

\bL{L:2.2}  Let $P$, $Q$ be rank-one projections in $B(H)$ and let  $a\le c\le b$. Then there are  projections $P'\sim Q'\sim P$  for which
$a P+b Q= c P'+(a+b-c)Q'$.
\eL
 \bp
The cases when $P=Q$ or when $a=0$ or $c=a$ or $c=b$ being  trivial, we assume that $P\neq Q$ and that  $0 < a<c<b$. Diagonalize the positive rank-two operator $A$. Then $A= a'E + b'F$ where $E$ and $F$ are two mutually orthogonal rank-one projections,  $0< a' \le a< c < b \le b'$, and $a'+b'= a+b$. Without loss of generality we can assume that $c=1$, and now  the conclusion follows from Lemma \ref {L:2.1}.

\ep
A generalization of Lemma \ref {L:2.2} provides the  algorithm for constructing frame perturbations in \cite {KvLd08}.

  The following lemma is obtained by iterative applications of  Lemma \ref {L:2.1}  and serves several
  complementary purposes: it illustrates in the simpler finite-dimensional case a
  construction that is applicable also in the cases of infinite dimensions,  it is  a
   key ingredient in the proof of Theorem \ref {T:4.3},  and it
  provides another  proof of Fillmore's characterization of finite sums of rank-one projections \cite [Theorem 1] {Fp69}.

\bL{L:2.3}  Let $A\in B(H)^+$ be a  finite rank operator and set 
\[
A=  \sum_{j=1}^n(1+\mu_j)(e_j\otimes e_j) +
\sum_{i=1}^m(1-\lambda_i)(f_i \otimes f_i)
\]
  where $\{e_1, \dots,
e_n, f_1, \dots, f_m\}$  is a collection of mutually orthogonal unit
vectors,
 $\mu_j > 0$ and $0 \le \lambda_i < 1$ for all $1\le i\le m$ and $1\le j\le
 n$. If all
  the eigenvalues of $A$ are greater than $ 1$,  set $m=0$, i.e., drop the sum involving the $\lambda_i$.  Similarly,  if all the eigenvalues of $A$ are less than or equal to $1$, set $n=0$. 
\item[(i)] Assume that $k:=\tr(A)- \tr(R_A) \in \mathbb N\cup\{0\}$. Then $A$ is the sum of
 $n+m + k$ rank-one projections. 
\item [(ii)]  Assume that $0\le \sum_{j=1}^n \mu_j-  \sum_{i=1}^m\lambda_i\le \max\{\mu_j\}$.
 Then there are $n+m$ rank-one projections  $P_1, P_2, \cdots, P_{m+n}$ for which
\be{e:5}
A= \sum_{h=1}^{m+n-1}P_h + \Big(1+ \sum_{j=1}^n\mu_j -  \sum_{i=1}^m\lambda_i\Big)P_{m+n}.
\ee
\eL
\bp
First notice that  $A_+= \sum_{j=1}^n\mu_j(e_j\otimes e_j) $ and  $A_-=\sum_{i=1}^m
\lambda_i(f_i\otimes f_i)$, hence by \eqref  {e:1}
\[
\tr(A)- \tr(R_A)=\tr(A_+)- \tr (A_-)  =\sum_{j=1}^n\mu_j - \sum_{i=1}^m\lambda_i .
\]

\item[(i)] To avoid triviality, assume that $A\ne 0$, and  in particular that $n\ne
0$. If $k>0$, let
\[
A_1:= \mu_1(e_1\otimes e_1) +\sum_{j=2}^n(1+\mu_j)(e_j\otimes e_j )+
\sum_{i=1}^m
 (1-\lambda_i)(f_i \otimes f_i).
\]
 Then $A= e_1\otimes e_1+A_1$,  $A_1\ge 0$, $R_{A_1}=R_A$, and $\tr(A_1)=\tr(A)-1$, whence\linebreak
 $\tr(A_1)-\tr(R_{A_1})=k-1. $  Iterating, we decompose $A$  into the sum of
 $k$ rank-one projections and a positive operator $A_k$ with $\tr(A_k)=\tr(R_{A_k})$.
  Thus, we can simply assume that $k=0$.  Hence  $ \sum_{j=1}^n\mu_j= \sum_{i=1}^m\lambda_i $ and also $m\ne0$.
  
Now we  start by the decomposition
\[
(1-\lambda_1)(f_1 \otimes f_1)+ (1+\mu_1) (e_1\otimes e_1)= P_1+ (1+\delta_1)(v_1\otimes v_1)
\]
where $\delta_1:= \mu_1-\lambda_1$ and $P_1$ and $v_1\otimes v_1$ are the rank-one
projections prescribed by Lemma \ref {L:2.1}. Then either  $\delta_1=0$ and $n=m=1$, in which case $A$ is the sum of two rank-one projections, or one of the three conditions hold: $\delta_1>0$, in which case $m > 1$;  $\delta_1<0$, in which case $n > 1$; or $\delta_1=0$ and $(n,m)\ne (1,1)$, in which case $n>1$ and $m >1$. Notice that  $v_1$ is  orthogonal to each $e_j$ and $f_i$ for $j>1$ and $i>1$ if any, so  Lemma \ref {L:2.1} yields again the decomposition
\[
\begin{cases}
(1+\delta_1)(v_1\otimes v_1)+(1-\lambda_2)(f_2 \otimes f_2) \quad &\text{if } \delta_1 > 0 \\
(1+\delta_1)(v_1\otimes v_1)+ (1+\mu_2)(e_2 \otimes e_2)   &\text{if } \delta_1 < 0\\
 (1-\lambda_2)(f_2 \otimes f_2) + (1+\mu_2)(e_2 \otimes e_2) &\text{if }  \delta_1 = 0
\end{cases}\quad = \quad  P_2+ (1+\delta_2)(v_2\otimes v_2) 
\]
where $P_2$ and $v_2\otimes v_2$ are rank one projections and 
\[
\delta_2= \begin{cases}
\mu_1-\lambda_1-\lambda_2 \quad  &\text{if } \delta_1 > 0 \\
 \mu_1+\mu_2 -\lambda_1 &\text{if } \delta_1 < 0\\
 \mu_2-\lambda_2&\text{if }  \delta_1 = 0
\end{cases}
\]

In general after $q$  steps, we have
\be{e:6}
A_q:= \sum_{j=1}^{n'}(1+\mu_j)(e_j\otimes e_j) + \sum_{i=1}^{m'} (1-\lambda_i)(f_i \otimes
f_i)=\sum_{j=1}^qP_j+ (1+\delta_q)(v_{q}\otimes v_{q})
\ee
 where $\delta_q=  \sum_{j=1}^{n'}\mu_j- \sum_{i=1}^{m'}\lambda_i$ and $n', m'\in \mathbb N$, $n'\le n, m'\le m$.
We continue the process until we ``run out"  of summands to which  apply Lemma \ref {L:2.1}. This occurs only when $n'=n$ and $m'=m$.  Indeed, if $n'=n$ but  $m' < m$, then $\delta_q= \sum_{j=1}^n\mu_j-\sum_{i=1}^{m'}\lambda_i= \sum_{i=m'+1}^m\lambda_i >0$ and thus we can further decompose
$(1+\delta_q)v_{q}\otimes v_{q}+(1-\lambda_{m'+1})(f_{m'+1}\otimes f_{m'+1})$ into the sum of a rank-one projection and a positive remainder. The case when $m'=m$ but $n'\ne n$ is similar. But when $n'=n$ and $m'=m$, then $\delta_q=0$ and hence $A= A_q$ is the sum of  $\tr(A)=n+m+k$ rank-one projections. 
\item[(ii)]  Assume without loss of generality that $\max\{\mu_j\}$ occurs for $j=n$, i.e.,
that
\[
 \sum_{j=1}^{n-1} \mu_j\le \sum_{i=1}^m\lambda_i\le  \sum_{j=1}^n \mu_j.
\]
We can carry on the same construction process as in  (i). If after  the $q$ steps that lead to the decomposition  \eqref{e:6}  we have $n'=n$ and $m'\ne m$, then 
 $\delta_q\ge \sum_{i=m'+1}^m\lambda_i \ge 0$ and we can continue the process. If we have   $m'=m$ but $n'\ne n$ then
\[
\delta_q=\sum_{j=1}^{n'}\mu_j-\sum_{i=1}^{m}\lambda_i\le \sum_{j=1}^{n-1}\mu_j-
\sum_{i=1}^{m}\lambda_1\le 0.
\]
and in this case too we can continue the process. Thus the process terminates only when $n'=n$ and $m' = m$ and thus (\ref {e:5}) holds.

\ep
\bR{R: 2.4} The condition in (ii) is necessary, because  if \eqref {e:5} holds, then 
\[
1+ \sum_{j=1}^n\mu_j -  \sum_{i=1}^m\lambda_i\le ||A||= 1+ \max \{\mu_j\}.
\]
\eR

This lemma provides  a constructive  proof  of  Fillmore's characterization of finite sums of finite projections \cite [Theorem 1] {Fp69}  that does not depend on the mean value theorem (see also \cite [Proposition 6] {DFKLOW}).
\bC{C:2.5}  \cite [Theorem 1] {Fp69} Let $A\in B(H)^+$ be a  finite rank operator.
Then $A$ is the sum of projections if and only if $\tr A \in \mathbb N$ and
 $\tr A \ge\tr (R_A).$
\eC
\bp
The sufficiency is given by Lemma \ref {L:2.3}  (i). For the necessity, assume that $A= \sum_{j=1}^k P_i$ is a sum of  projections and by further decomposing them if necessary, assume that they all have rank one. Then $\tr A=k\in\mathbb N $ and, clearly,  $\text{rank }A\le k$. .

\ep

The matrix construction in Lemma \ref {L:2.1}  extends to $C^*$-algebras and hence in particular to von Neumann algebras. It is well known that  given a  collection
$\{E_j\}_{j=1}^n$ of mutually orthogonal equivalent projections in
a $C^*$-algebra $\mathscr A$, we can chose a corresponding set of matrix units and hence an embedding of $M_n(\mathbb C)$ into $\mathscr A$. Thus by Lemma
\ref {L:2.1} and Corollary \ref{C:2.5}  we obtain:

\bL{L:2.6}  Let $\mathscr A$ be a $C^*$-algebra.
\item[(i)]
If $E$ and $F$ are two mutually orthogonal equivalent projections in $\mathscr A$,
 $0\le\lambda\le 1$, and $\mu \ge 0$, then there are two projections $P_-$ and $P_+$
 in  $\mathscr A$, with $P_-\sim P_+\sim E$,  for which $(1+\mu) E+(1-\lambda)
 F= P_-+(1+\mu-\lambda)P_+.$
\item [(ii)] If $A=\sum_{j=1}^n\gamma_j E_j$ for some mutually orthogonal equivalent
projections $E_j\in \mathscr A$ and some scalars $\gamma_j > 0$ with
$\sum_{j=1}^n\gamma_j =k\in \mathbb N$ and $k\ge n$, then $A$ is the sum of $k$
equivalent projections in $\mathscr A$.
\eL

\bR{R:2.7} 
\item[(i)] The embedding  of $M_n(\mathbb C)$ into $\mathscr A$ depends not only
on the projections $E_j$ but also on the  matrix units.
However, once these
 matrix units are chosen, the construction in Lemma \ref {L:2.1}  assigns the
 decomposition in a canonical way.

Explicitly for the $n=2$ case, let $V\in \mathscr A$ be a partial isometry with
 $E=V^*V$ and $F=VV^*$, then the projections $P_-$ and $P_+$ obtained from this
 embedding and the formulas in Lemma  \ref {L:2.1} are
\begin{equation}\label {e:7}
\begin{split}
P_-&:=(1-\rho) E-\sqrt{\rho(1-\rho)}(V+V^*)+\rho F  \\
P_+&:=(1-\nu) E+\sqrt{\nu(1-\nu)}(V+V^*)+\nu F.
\end{split}
\end{equation}
It can also be verified directly that setting $W:=\sqrt{1-\rho}E-\sqrt{\rho}V$, we get  $W^*W=E$ and $WW^*=P_-$. Thus $W$ is a partial isometry and hence $P_-$ is a projection and
 $P_-\sim E$. Similarly, $P_+\sim E$.

\item[(ii)] More generally, if  $0\le a\le c\le b$,  set
\[
\rho:=\begin{cases}\frac{a(b-c)}{c(b-a)}\quad&\text{if } b\ne  c\\0&\text{if } b=c \end{cases}
\quad \quad \text{and}\quad \quad
 \nu:= \begin{cases}\frac{a(c-a)}{(a+b-c)(b-a)}\quad&\text{if } b\ne  c\\1&\text{if } b=c\end{cases}.
\]
Then  with $P_-$ and $P_+$ as in (\ref{e:7}), $bE+aF= cP_-+(a+b-c)P_+$.
\eR

\bL{L:2.8}  Let $P$, $Q$ be finite equivalent commuting projections in a  von Neumann algebra $M$ and let $0\le a < b$  and $a\le c \le b$. Then there are  projections $P'\sim Q'\sim P$ in  $M$ for which
$a P+b Q= c P'+(a+b-c)Q'$.
\eL
\bp
By the assumption of finiteness, we have the cancellation $P-PQ\sim Q-PQ$. By Lemma \ref {L:2.6}   (see also  Remark \ref {R:2.7}  (ii)), there are projections $P'\sim Q'\sim P-PQ$ in $M$, with $P'\vee Q'\le P-PQ+ Q-PQ$, for which $a (P-PQ)+b (Q-PQ)= c P'+(a+b-c)Q'$. But then,
\[
a P+bQ= a (P-PQ)+b (Q-PQ) +  (a+b)PQ
= c(P'+PQ)+ (a+b-c)(Q'+PQ).
\]
Since $P'\perp PQ$ and $Q'\perp PQ$, $P'+PQ$ and $Q'+PQ$ are projections and both are equivalent to $P$.

\ep
\bR{R:2.9} 
\item[(i)] If the projections $P$ and $Q$ are not finite, cancellation might fail and indeed the property itself might fail. For instance if $P$ is infinite but $P\ne I$, then there are no projections $P'$ and $Q'$ for which $\frac{1}{5}P+ I =  \frac{2}{5}P'+  \frac{4}{5}Q'$. Indeed, otherwise
$I-Q'=\frac{2}{5}P' - \frac{1}{5}P- \frac{1}{5}Q'$, whence $||I-Q'||\le \frac{4}{5}$ and hence $Q'=I$. But then, $\frac{1}{5}P+ \frac{1}{5}I= \frac{2}{5}P'$, whence $P=P'=I$, against the assumption.
\item[(ii)] Lemma \ref {L:2.8}   holds also for every $C^*$-algebra $\mathscr A$ with the cancellation property (e.g., AF-algebras.)
\eR

The proof of Lemma \ref {L:2.1}  can be used also to obtain a simple constructive proof of Fillmore's \cite [Theorem 2, Corollary] {Fp69} characteriziation of  the operators in $B(H)$ that are sums of two projections. The same characterization holds for von Neumann algebras.

\bP {P:2.10}  Let $M$ be a von Neumann algebra and $A\in M$ with $0\le A\le 2I$. Then $A$ is the sum of two projections in $ M$ if and only  $A=E\oplus B$ where $E$ is a (possibly zero) projection in $M$ and there is a unitary $U\in M$ that commutes with $E$ and for which $UBU^*=2R_B-B.$

 \eP
 \bp
To prove the sufficiency, it is obviously enough to consider the
case when $E=0$ and $R_A=I$, i.e., $UAU^*=2I-A$. Let $\mathscr A$ be
the (abelian) von Neumann algebra generated by $A$. Since $UAU^*=2I-A\in
\mathscr A$, it follows that $U\mathscr A U^*\subset \mathscr A$ and hence $U^2 \mathscr A(U^2)^*\subset U\mathscr A U^*$. Since
\[
A=2I- UAU^*= 2I-U(2I-UAU^*)U^*= U^2A(U^2)^*,
\]
it follows that $ \mathscr A = U^2 \mathscr A(U^2)^*$.  Thus
 $\mathscr A=U\mathscr A U^*$, i.e., $U\cdot U^*$ is a
conjugation of    $\mathscr A$. In particular, for every Borel set
$\Omega\subset [0,2]$ there is a Borel set $\Omega_U\subset [0,2]$
for which \linebreak $U\chi_A(\Omega)U^*= \chi_A(\Omega_U)$.

Let $E_t:=\chi_A[0, t)\in \mathscr A$   for $t\in [0,2]$  be the spectral resolution of $A$.  Then
\ba
A&=  \int_{[0,1)} tdE_t + \chi_A\{1\} + \int_{(1,2]} tdE_t = U(2I-A)U^* \\
&=U\Big( \int_{(1,2]} (2-t)dE_t \Big)U^* +U \chi_A\{1\} U^*+ U\Big(
\int_{[0,1)} (2-t)dE_t  \Big)U^*.
\end{align*}
It is now clear that 
\[
U\chi_A\{1\}U^*= \chi_A\{1\}, \quad U\chi_A[0,1)U^*=
\chi_A(1,2],\quad \text{and} \quad U\chi_A(1,2]U^*= \chi_A[0,1).
\]
 In particular, $
\int_{(1,2]} tdE_t = U\Big(  \int_{[0,1)} (2-t)dE_t\Big )U^*. $ Thus
\[
A =  \int_{[0,1)} tdE_t +\chi_A\{1\}+ U\big( \int_{[0,1)} (2-t)dE_t\big)U^*.
\]
Let
\ba
2P_-:&= \int_{[0,1)} tdE_t-U \int_{[0,1)}\sqrt{t(2-t)}dE_t - \int_{[0,1)} \sqrt{t(2-t)}dE_t U^*+ U \int_{[0,1)}(2-t)dE_tU^*\\
2P_+:&=\int_{[0,1)}tdE_t+U \int_{[0,1)}\sqrt{t(2-t)}dE_t +  \int_{[0,1)}\sqrt{t(2-t)}dE_t U^*+ U \int_{[0,1)}(2-t)dE_tU^*.
\end{align*}

Then both $P_-$ and $ P_+$ belong to $M$ and are selfadjoint.  Since
$\chi_A[0,1)\perp U\chi_A[0,1)U^*$, it is simple to verify that
$P_-, P_+$ are idempotents and hence are projections. Furthermore $
P_-+P_+=A-\chi_A\{1\}$, hence $P_-\perp \chi_A\{1\}$ and thus
$P_-+\chi_A\{1\}$ is also a projection, which completes the proof of
the sufficiency. The necessity follows as in Fillmore's proof in
\cite [Theorem 2, Corollary] {Fp69} from  the analysis of the
relative position of two projections which holds for general von
Neumann algebras (e.g., see \cite [Pgs 306-308]{TOAI}), and hence,
applies without changes to our setting.
 \ep

\bR{R:2.11}  With the notations of the above proof, if  $ \mathscr
A$ is a masa,  then it cannot be singular, since $U$ belongs to the
normalizer $\mathscr N(\mathscr A)$ of $ \mathscr A$ but does not
belong to $ \mathscr A$, as otherwise $A=I$, against the assumption
that $ \mathscr A$ is a masa. \eR

\section{The necessary condition}

\bP{P3.1} 
Let $A\in M^+$ and let $N\in \mathbb N\cup\{\infty\}$. 
Then the following conditions are equivalent.
\item[(i)] There is a partial isometry $V$ with $V^*V= R_A$ and  a decomposition of the identity   into $N$ mutually orthogonal nonzero projections $E_j$, $I=\sum_{j=1}^NE_j$,   for which $\sum_{j=1}^NE_jVAV^*E_j= I$,  the convergence of the series being in the strong operator topology if $N=\infty$.
\item[(ii)] $A$ is the sum of $ N $ nonzero projections, the convergence of the series being in the strong operator topology if $N=\infty$, and if $M$ is semifinite, then $\tau(A)=\tau(I)$.
\eP

\bp
\item[(i)] $\Longrightarrow$ (ii)  For every $j$, let  $W_j:=E_jVA^{\2}$ and let $P_j:=W_j^*W_j$. Since
\[
W_jW_j^*= E_jVAV^*E_j= E_j\sum_{i=1}^NE_iVAV^*E_iE_j = E_j,
\]
we see that $W_j$ is a partial isometry, and hence, $P_j$ is a
projection and $P_j\sim E_j$ for every $j$. Then
\[
\sum_{j=1}^N P_j = \sum_{j=1}^N A^{\2}V^*E_jVA^{\2} = A^{\2}V^*\big(\sum_{j=1}^N E_j\big) VA^{\2}= A^{\2}V^*VA^{\2}=A^{\2}R_AA^{\2}= A
\]
and if $N=\infty$ the series $\sum_{j=1}^N E_j$ and hence the series $\sum_{j=1}^N P_j$ converge in the strong operator topology.
Furthermore, if $M$ is semifinite, by the normality of the trace $\tau$ we have
\[
\tau(A)= \sum_{j=1}^N\tau(P_j)= \sum_{j=1}^N\tau(E_j)= \tau(I).
\]

\item[(ii)] $\Longrightarrow$ (i)  Let $A= \sum_{j=1}^NP_j$ where $\{P_j\}$ are nonzero projections. First, we decompose the identity $I=\sum_{j=1}^NE_j$ into $N$ mutually
orthogonal projections $E_j\sim P_j$. This is immediate if all the
projections $P_j$ are infinite, and hence  so is $I$, because then
we can decompose $I$ into  $N$ mutually orthogonal infinite
projections and all infinite projections are equivalent by the assumption that $M$ is $\sigma$-finite.  Assume
henceforth that $M$ is semifinite and that $\Lambda:=\{j\mid \tau(P_j) <
\infty\}\ne \emptyset$ and let $\Lambda'$ be its  (possibly empty) complement.
Then
\[
 \tau(I)=\tau(A) = \sum_{j=1}^N\tau(P_j)\ge \sum_{j\in \Lambda}\tau(P_j).
 \]
Whether $M$ is of type II or it is of type I and then $ \sum_{j\in \Lambda}\tau(P_j)\in\mathbb N \cup \{\infty\}$, there exists a projection $F$  with $\tau(F)= \sum_{j\in \Lambda}\tau(P_j)$. Then it is routine to find mutually orthogonal projections $E_j\le F$ with $\tau(E_j)=\tau(P_j)$ for every $j\in \Lambda$. Let $E:=  \sum_{j\in \Lambda} E_j$. Then $E\le F$ and $\tau(E)=\tau(F)=\tau(I)$.   We now consider three cases.

In the first case,  assume that $\tau(I)< \infty$. Then $\tau(I-E)=0$, hence $E=I$, and we are done.

In the second case, assume that  $\tau(I)= \infty$ and $\Lambda' = \emptyset$. Then $E$ is an infinite projection, hence there is an isometry $W$ for which $WW^*=E$. Set $E'_j=W^*E_jW$. Then $ E'_j\sim E_j\sim P_j$ for every $j$ and $I= \sum_{j=1}^NE'_j $ provides the required decomposition.

In the third case, assume that $\tau(I)= \infty$ and $\Lambda' \ne \emptyset$. Modify if necessary $F$ so that $I-F$ is infinite  and hence  so is  $I-E\ge I-F$. Then  decompose $I-E$
into  $\card\, \Lambda'$ mutually orthogonal infinite projections $E_j$, $I-E= \sum_{j\in \Lambda'}E_j$. Since  $P_j\sim E_j$ for  all $ j\in \Lambda'$, $I= \sum_{j\in \Lambda} E_j+ \sum_{j\in \Lambda'} E_j=\sum_{j=1}^NE_j $
 provides in this case too the required decomposition.

Now choose  partial isometries $W_j$ with $P_j=W_j^*W_j$ and
$E_j=W_jW_j^*$.  If $N<\infty$, set $B:=  \sum_{j=1}^N W_j$.  If
$N=\infty$ and $m> n$, then 
\be{e:8}
\big(\sum_{j=m}^nW_j\big)^*\big(\sum_{j=m}^nW_j\big)=
\sum_{i,j=m}^nV_i^*W_j= \sum_{j=m}^nW_j^*W_j = \sum_{j=m}^nP_j. \ee
Thus, by the strong (and hence the weak) convergence of the series
$\sum_{j=1}^\infty P_j $, we see that the series $\sum_{j=1}^\infty
W_j $ is strongly Cauchy  and hence  converges in the strong
operator topology.  Again, call  its sum $B$. By the same
computation as in \eqref {e:8}, we have $B^*B=\sum_{j=1}^N P_j =A.$
 Let $B=VA^{\frac{1}{2}}$ be the polar decomposition of $B$.
Then $V^*V= R_A$ and $BB^*=VAV^*$. Moreover, $E_jB=W_j$ for every
$j$, thus
\[
\sum_{j=1}^NE_jVAV^*E_j = \sum_{j=1}^NE_jBB^*E_j= \sum_{j=1}^NW_jW_j^* = \sum_{j=1}^NE_j = I.
\]
\ep

\bL{L3.2} Let $A\in M^+$ be a strong sum of projections.
\item [(i)] Either $||A||>1$ (equivalently,  $A_+\ne 0$) or $A$ is a projection.
\item [(ii)] If $M$ is semifinite, then $\tau(R_A)\le \tau(A)$.
\eL
\bp
\item [(i)] Obvious, since if $P$, $Q$ are projections, then $||P+Q||= 1$ if and only if $PQ=0$ if and only if $P+Q$ is a projection.
\item [(ii)] Let $A= \sum_{j=1}^N P_j$ with $P_j$ nonzero projections and $N\in \mathbb N\cup\{\infty\}$ and assume without loss of generality that $ \tau(A)< \infty$. By Kaplansky's parallelogram law, \cite[Theorem, 6.1.7] {KrRj2}, for every integer $n\le N$ we have
\[
\tau\big( \bigvee_{j=1}^nP_j\big)\le \sum_{j=1}^n\tau(P_j)\le
\tau(A). \]
If $N< \infty$, then $R_A= \bigvee_{j=1}^NP_j$ and we are done. If $N= \infty$, then
$\bigvee_{j=1}^nP_j \uparrow R_A$ and by the normality of $\tau$,  $\tau\big(\bigvee_{j=1}^nP_j\big)\uparrow \tau(R_A)$. Thus also $ \tau(R_A)\le \tau(A)$.

\ep

\bT{T:3.3} 
 Assume that $A \in M$ is a strong sum of projections. Then
 \item[(i)]  If $M$ is  of type I, then  $\tr(A_+)\ge\tr(A_-)$ and  either  $\tr(A_+)=\infty$ or $ \tr(A_+)<\infty$  and  $\tr(A_+)-\tr(A_-)\in \mathbb N \cup\{0\}$.
 \item[(ii)] If $M$ is  of type II,  then $\tau(A_+)\ge \tau(A_-)$.
 \item[(iii)] If  $M$ is of type III, then either $||A||>1$ (equivalently,  $A_+\ne 0$) or $A$ is a projection.
\eT

\bp 
 \item[(iii)] is given by Lemma  \ref {L3.2}  (i), so assume henceforth that $M$ is semifinite.
Let $A= \sum_{j=1}^N P_j$ with $P_j$ nonzero projections and $N\in \mathbb N\cup\{\infty\}$.

Assume first that  $\tau(R_A) < \infty$ and hence also  $\tau(A)< \infty$ and $\tau(A_-)< \infty$. Then by (\ref {e:1}) and by Lemma \ref {L3.2} we have 
$\tau(A_+)- \tau(A_-) = \tau(A) - \tau(R_A) \ge 0.$
Moreover, if $M$ is of type I, then $N < \infty$ and both $\tau(A)$ and $\tau(R_A)$ are positive integers, which proves the integrality condition in (i) for the case when $\tau(R_A) < \infty$ (see  also  Corollary \ref{C:2.5} ).

Now assume that $\tau(R_A) =\infty$  and assume furthermore that  $\tau(A_+) < ~\infty$.  Obviously,  $\tau(I) =\infty$  and by Lemma \ref{L3.2} , $\tau(A) = \infty$, hence $\tau(A) =\tau(I)$. Thus  by   Proposition \ref {P3.1} there is  a partial isometry $V$ with $V^*V= R_A$ and  a
decomposition of the identity $I=\sum_{j=1}^NE_j$  into $N$ mutually
orthogonal projections $E_j$  for which $\sum_{j=1}^NE_jVAV^*E_j=
I$.  Recall that the map \[M\ni X\to \Phi(X):= \sum_{j=1}^N E_jXE_j
\in M\] is linear, positive, unital, faithful, and in case $M$ is
semifinite, it is also trace preserving. Then we have by (\ref
{e:1}) that
\[
I= \Phi(VAV^*)=\Phi(VA_+V^*)- \Phi(VA_-V^*)+\Phi(VR_AV^*).
\]
It follows from $V^*V=R_A$ that 
\be{e:9}
 \Phi(VA_+V^*)=\Phi(VA_-V^*)+ \Phi(I-VV^*). 
\ee
and hence 
\[
\tau(VA_-V^*)= \tau(\Phi(VA_-V^*))\le \tau(\Phi(VA_+V^*))= \tau(VA_+V^*)  =\tau(R_AA_+)=\tau(A_+)< \infty.
\]
But then,
\[
 \tau(A_-)=\tau (R_AA_-R_A)= \tau (V^*VA_-V^*V)= \tau (VV^*VA_-V^*)= \tau(VA_-V^*).
\]
This concludes the proof of the case when $M$ is of type II.
If $M$ is of type I and $\tr(A_+)< \infty$, it follows  from
 \eqref {e:9} and  the above computations that
\[
\tr(A_+)= \tr(A_-)+\tr(\Phi(I-VV^*))= \tr(A_-)+\tr(I-VV^*).
\]
This shows that $\tr(I-VV^*)< \infty$, i.e., $I-VV^*$ is a finite
projection, and therefore  $\tr(I-VV^*)\in \mathbb N \cup\{0\}.$

\ep

 \bR{R:conditions} Notice that if $M$ is semifinite and $A\in M^+$, then
 \[
 \tau(A_+)\ge \tau(A_-)\Longrightarrow \tau(A) \ge \tau(R_A)\not \Longrightarrow \tau(A_+)\ge \tau(A_-). \]
 However, if $\tau(R_A)<\infty$, then $\tau(A_+)\ge \tau(A_-)\Longleftrightarrow \tau(A) \ge \tau(R_A)$.
 \eR

\bR{R:3.4} In the case  of $M=B(H)$,  let $\mathscr {V}(A):=
\{VAV^*\mid V^*V=R_A\}$ be the partial isometry orbit of $A$ and let
$E$ denote the (unique) normal conditional expectation on the
diagonal masa of $B(H)$ (according to a fixed orthonormal basis).
Then Proposition \ref {P3.1}  states that a positive operator $A\in B(H)$
with infinite trace is a strong sum of rank-one projections if and
only if $I\in E(\mathscr {V}(A))$.  When  $A$ is also invertible,
this is a special case of \linebreak \cite [Proposition 4.5] {AMRS}.  
In the case of compact operators,  the diagonals of the partial isometry orbit  are characterized in terms of majorization of sequences by the infinite dimensional Schur-Horn theorem obtained  in \cite
{KW08}. The set $E(\mathscr {V}(A))$ is further studied in in \cite {KvLd09} for the  case of positive not necessarily compact operators. \eR

In the case of $M=B(H)$, an application of  Proposition \ref {P3.1}  together with a modification of the proof of Theorem \ref {T:3.3}  (i) provides an alternative proof of the necessity of Kadison's integrality condition in  \cite [Theorem 15]{Kr02b} that characterizes the diagonals of infinite co-infinite projections and identifies explicitly the integer as the difference of the traces of two projections.

\bC{C:3.6} \cite [Theorem 15]{Kr02b}
Let $P\in B(H)$ be an infinite, co-infinite projection, let $e_n$ be an orthonormal basis, let $c_n:= (Pe_n,e_n)$,  and assume that $\sum\{c_n \mid c_n \le \2\} <\infty$ and $\sum\{1-c_n \mid c_n > \2\} < \infty$.  Then  
\[
 \sum\{1-c_n \mid c_n > \2\}- \sum\{c_n \mid c_n \le \2\} \in \mathbb Z.
  \]
  \eC
\bp
Let $W$ be an isometry with $P=WW^*$. Define 
\[
w_n= \begin{cases} \frac{1}{\sqrt {c_n}}W^*e_n \quad &\text{if } c_n\ne 0\\
e_1&\text{if } c_n= 0
\end{cases}\quad\quad \text {and} \quad\quad  P_n:= w_n\otimes w_n
\]
Then $||w_n||=1$ for every $n$ and hence $P_n$ are rank-one projections. A simple computation shows that  $I= \sum_n c_nP_n$, with the series converging in the strong operator topology.
Define
\[
T_+:=\sum \{(1-c_n)P_n\mid c_n > \2\}, \quad\quad T_-:=\sum \{c_nP_n\mid c_n \le \2\}, \quad\quad {and} \quad\quad  T:=T_+-T_-.
\]
Then $T_+, T_-$, and hence $T$ are trace class operators and 
\[
\tr T =  \sum\{1-c_n \mid c_n > \2\}- \sum\{c_n \mid c_n \le \2\}. 
\]
Since $\sum \{ c_nP_n\mid c_n > \2\}$ and  $\sum \{ (1-c_n)P_n\mid c_n > \2\}$ both converge in the strong operator topology,  it follows that also $\sum \{ P_n\mid c_n > \2\}$  converges in the strong operator topology. Set $A:=\sum \{ P_n\mid c_n > \2\}$. By Proposition  \ref {P3.1}   there is a partial isometry  $V$ with $V^*V=R_A$ for which $E(VAV^*)=I$, where $E$ is the conditional expectation on the atomic masa (the operation of taking the main diagonal).

Since  $I= A-T$, we have $VV^*= VAV^*-VTV^*$ and thus
\[
E(VV^*)= E(VAV^*)-E(VTV^*)=I-E(VTV^*)= E(I)-E(VTV^*).
\]
Thus $ E(VTV^*)= E(I-VV^*)$ and hence 
\[
\tr(V^*VT) = \tr(VTV^*)= \tr(I-VV^*)\in \mathbb Z
\]
 since $T$ is trace-class and $VV^*$ and hence $I-VV^*$ are projections. On the other hand,  
 \[
 (V^*V)^\perp T = (V^*V)^\perp (A-I) = -(V^*V)^\perp,
 \]
  hence $\tr ((V^*V)^\perp T) \in \mathbb Z$, and thus 
\[
\tr (T) = \tr (V^*V T) + \tr ((V^*V)^\perp T) \in \mathbb Z.
\] 
\ep

\section{B(H): the finite trace case}
In this section we will prove that if   $A\in B(H)^+$ and  $\tr(A_-)\le 
\tr(A_+) < \infty$ and $\tr(A_+) -
\tr(A_-) \in \mathbb N \cup \{0\}$, then $A$ is a strong sum of
projections.  Since the trace-class operators $A_+$ and $A_-$ are
diagonalizable and have orthogonal supports, then by \eqref
{e:1}, $A$ too is diagonalizable. As in Lemma \ref {L:2.3} ,  let us denote the
eigenvalues of $A$ which are larger than $1$, if any, by $1+\mu_j$
and those are less or equal than $1$, if any, by $1-\lambda_i$,
i.e., set
\[
A:=  \sum_{j=1}^N(1+\mu_j)(e_j\otimes e_j) + \sum_{i=1}^K (1-\lambda_i)(f_i \otimes f_i)
\]
where $N, K\in \mathbb N\cup \{0\}\cup\{\infty\}$, the unit vectors $e_j, f_i$ are 
mutually orthogonal, $\mu_j > 0$, and $0
\le \lambda_i \le 1$ for all $i$ and $j$. Notice that the series, if infinite, converge
in the strong operator topology.  Of course, it would
 be equivalent to assume that $\mu_j\ge 0$ and $0< \lambda_i< 1$ for all $i$ and $j$.
  Thus $A_+= \sum_{j=1}^N\mu_j(e_j\otimes e_j) $ and $A_-=\sum_{i=1}^K \lambda_i(f_i \otimes f_i)$
  and hence \[ \sum_{j=1}^N\mu_j- \sum_{i=1}^K \lambda_i\in \mathbb N \cup \{0\}.\]
Here too we adopt the convention to set a series $\sum_{i=1}^0$ as
zero, e.g., by $K=0$ we
 mean that $A$ has no non-negative eigenvalues less or equal than $1$, and hence, $A_-=0$;
 similarly for $N=0$.

Our proof will depend  on iterative applications of Lemma \ref {L:2.1}. Since we will focus on
infinite rank operators, i.e., on the case when $N+K=\infty$, the process will not terminate as in Lemma
 \ref {L:2.3}  after a finite number of steps and the crux of the proofs will be to establish strong convergence.
This will be illustrated by the following lemma which handles two key special cases.

\bL{L:4.1} Let $\{g_o, g_1, \dots \}$ be mutually orthogonal unit vectors.
\item [(i)] Let $A= (1-\lambda) (g_o \otimes g_o) + \sum_{j=1}^\infty (1+\mu_j) (g_j \otimes g_j)$
where $ \mu_j >0$,  $0 \le \lambda \le 1$,  for all $j$ and $\lambda = \sum_{j=1}^\infty \mu_j$.
Then $A$ is a strong sum of projections.
\item[(ii)] Let $A= (1+\mu) (g_o \otimes g_o) + \sum_{j=1}^\infty (1-\lambda_j) (g_j \otimes g_j)$
where  $ \mu > 0 $, $ 0\le \lambda_j \le 1$ for all $j$ and $\mu = \sum_{j=1}^\infty \lambda_j$.
Then $A$ is a strong sum of projections.
\eL

\bp
\item[(i)]  If  $\lambda=0$, then $\mu_j=0$ for all $j$ and hence $A$ is already a projection. Thus assume that $\lambda \ne 0$. Define
\[
\delta_j:=
\begin{cases} -\lambda\quad & j=1\\
\sum_{i=1}^{j-1}\mu_i -\lambda &j>1.
\end{cases}
\]
Then $\delta_j$ increases strictly to $0$, so we can also define
 \be {e:10}
\sigma_j:= \begin{cases} 0\quad& j=1\\ \frac{(1+\delta_{j-1})\delta_{j-1}}{(1+\delta_j)(2\delta_{j-1}-\delta_j)}& j > 1.
\end{cases}
\ee
Then for every $j> 1$,  $\sigma_j>0$   and also
\be{e:11}
1-\sigma_j=\frac{(1+\delta_j-\delta_{j-1})(\delta_{j-1}-\delta_j)}{(1+\delta_j)(2\delta_{j-1}-\delta_j)}> 0.
\ee
 Define also
\be{e:12}
v_j:=\begin{cases}g_o\quad&  j=1\\
 \sqrt{\sigma_{j}}v_{j-1}+  \sqrt{1-\sigma_{j}}g_{j-1} & j>1.
 \end{cases}
\ee
Solving this recurrence relation, we get
\be{e:13}
v_j= \sum_{k=0}^{j-2}\Big(\sqrt{1-\sigma_{k+1} }\prod_{i=k+2}^j \sqrt{\sigma_i}\Big)g_k + \sqrt{1-\sigma_{j}}g_{j-1}.
\ee
 We claim that there is a sequence of rank-one projections $P_j$ for which
 \be{e:14}
 (1-\lambda) (g_o \otimes g_o) + \sum_{j=1}^n(1+\mu_j) (g_j \otimes g_j) =  \sum_{j=1}^n P_j+ (1+\delta_{n+1})(v_{n+1}\otimes v_{n+1})
 \ee
 for every $n$.  By Lemma \ref {L:2.1} 
\[
(1-\lambda)(g_o \otimes g_o)  + (1+\mu_1)(g_1\otimes g_1) =P_1+(1+\mu_1-\lambda)(v\otimes v)= P_1+(1+\delta_2)(v\otimes v)
\]
where $P_1$ is a rank-one projection and by \eqref {e:3} and \eqref {e:2}, $v= \sqrt{\nu}g_o+\sqrt{1-\nu}g_1$ and
\[
\nu=  \frac{(1-\lambda)\lambda}{(1+\mu_1-\lambda)(\mu_1+\lambda)} = \frac{(1+\delta_1)(-\delta_1)}{(1+\delta_2)(\delta_2-2\delta_1)}= \sigma_2.
\]
Thus $v= v_2$ and hence (\ref {e:14}) is satisfied for $n=1$.
Assume that  (\ref {e:14}) is satisfied for $n-1$. Then

\begin{alignat*}{2}
(1-\lambda)&(g_o \otimes g_o)  + \sum_{j=1}^n(1+\mu_j)(g_j\otimes g_j)= \\
&=\sum_{j=1}^{n-1} P_j + (1+\delta_n)(v_n\otimes v_n)+ (1+\mu_n)(g_n\otimes g_n) \quad &(\text{by  the induction hypothesis})\\
&=\sum_{j=1}^{n-1} P_j + P_n + (1+\mu_n+\delta_n)(v\otimes v) & (\text{by Lemma \ref {L:2.1} })\\
&= \sum_{j=1}^{n} P_j  + (1+\delta_{n+1})(v\otimes v)  & (\text{by the definition of $\delta$} )
\end{alignat*}

where $P_n$ is a rank-one projection, and by \eqref {e:3} and \eqref{e:2}, $ v = \sqrt{\nu} v_n + \sqrt{1-\nu}g_n$ and
\[
\nu=  \frac{(1+\delta_n))(-\delta_n)}{(1+\mu_n+ \delta_n)(\mu_n-\delta_n)} = \frac{(1+\delta_n)\delta_n}{(1+\delta_{n+1})(2\delta_n-\delta_{n+1})}= \sigma_{n+1}.
\]
Hence $v=v_{n+1}$ and thus (\ref {e:14}) is satisfied for $n$. Thus for every $n$,
\[
A=  \sum_{i=1}^nP_i +  (1+\delta_{n+1})(v_{n+1}\otimes v_{n+1}) +  \sum_{j=n+1}^\infty (1+\mu_j) (g_j \otimes g_j).
\]
Since $\sum_{j=n+1}^\infty (1+\mu_j) (g_j \otimes g_j)\underset
{s}\to 0$, to prove that $A=  \sum_{i=1}^\infty P_i $ where the convergence is in
the strong topology (and hence,  to establish the thesis), we need
to show that $v_{n+1}\otimes v_{n+1}\underset {s}\to 0$, or,
equivalently, that $v_j\to 0$ weakly. Since $v_j\in
\text{span}\{g_i\}$, it is enough to show that $(v_j, g_q) \to 0$
for every  $q\in \mathbb N\cup \{0\}$. Indeed, for  every $j > q+1$,
we have from (\ref {e:13}) that
\[
(v_j, g_q) =  \sqrt{1-\sigma_{q+1} } \prod_{i=q+2}^j \sqrt{\sigma_i}.
\]
Thus it is enough to show that $\prod_{i=2}^j \sigma_i \to  0$, or, equivalently, that
 $\sum_{i=2}^\infty  (1-\sigma_i) = \infty$. By (\ref {e:11}) and since
 $\delta_{j-1}< \delta_j <0$ we have
\be{e:15}
1-\sigma_j=\frac{(1-\delta_{j-1}+\delta_j)(\delta_{j-1}-\delta_j)}{(1+\delta_j)(2\delta_{j-1}-\delta_j)}
> \frac{\delta_{j-1}-\delta_j}{2\delta_{j-1}-\delta_j } >  \frac{1}{2}\frac{\delta_{j-1}-\delta_j}{\delta_{j-1}}> 0.
\ee Since  $\delta_j\uparrow 0$, for every $n > m$, 
\[
\sum_{i=m+1}^n \frac{\delta_{i-1}-\delta_i}{\delta_{i-1}}
\ge\sum_{i=m+1}^n \frac{\delta_{i-1}-\delta_i}{\delta_m} =
\frac{\delta_{m} - \delta_n}{\delta_m}= 1-\frac{\delta_n}{\delta_m},
\]
 whence $\sum_{i=m+1}^\infty
\frac{\delta_{i-1}-\delta_i}{\delta_{i-1}}> \2$ for every $m$. As a
consequence, $\sum_{j=2}^\infty\frac{
\delta_{i-1}-\delta_i}{\delta_{i-1}}=\infty$, and thus,
$\sum_{j=2}^\infty(1-\sigma_j)=\infty$, which  completes the proof for
this case.

\item[(ii)]
Let $k:=\card \{j \mid \lambda_j=1\}$.  By passing to
\begin{align*}
A':&=A-k(g_o\otimes g_o)  -\sum\{ g_j \otimes g_j  \mid \lambda_j=0\}\\
&= (1+\mu-k) (g_o\otimes g_o)+\sum\{ (1-\lambda_j)(g_j \otimes g_j ) \mid  0<\lambda_j< 1\},
\end{align*}
we can assume without loss of generality that $0< \lambda_j<  1$ for all $j$.
Define
\[
\delta_j:= \begin{cases} \mu & j=1\\
\mu-\sum_{i=1}^{j-1}\lambda_i &j>1
\end{cases}.
\]
Then $\delta_j \downarrow 0$. Let $\sigma_j$ and $v_j$ be defined by (\ref {e:10})
 and (\ref {e:12}) respectively.
 We claim that there is a sequence of rank-one projections $P_j$ for which

 \be{e:16}
 (1+\mu) (g_o \otimes g_o) + \sum_{j=1}^n(1-\lambda_j) (g_j \otimes g_j) =
  \sum_{j=1}^n P_j+ (1+\delta_{n+1})(v_{n+1}\otimes v_{n+1})
 \ee
 for every $n$.

 Apply Lemma \ref {L:4.2} to obtain
 \[
 (1-\lambda_1) (g_1\otimes g_1)+ (1+\mu) (g_o\otimes g_o) = P_1+ (1+\delta_2)(v\otimes v )
 \]
where $P_1$ is a rank-one projection and by  \eqref {e:3}, \eqref {e:2}, and \eqref {e:11},
$v = \sqrt{\nu} g_1 + \sqrt{1-\nu } g_o$ and
\[
\nu=  \frac{(1-\lambda_1)\lambda_1}{(1+\mu-\lambda_1)(\mu+\lambda_1)} =
\frac{(1+\delta_2-\delta_1)(\delta_1-\delta_2)}{(1+\delta_2)(2\delta_1-\delta_2)}= 1-\sigma_2.
\]
Thus $v=v_2$ and (\ref {e:16}) holds for $n=1$. The
inductive proof of the claim then proceeds as  in part (i). Thus
\[
A = \sum_{j=1}^nP_j+ (1+\delta_{n+1})(v_{n+1}\otimes v_{n+1}) +
\sum_{j=n+1}^\infty(1-\lambda_j)(g_j\otimes g_j),
\]
and hence, to prove that $A=\sum_{j=1}^\infty P_j$ we need to show
that $v_j\to0$ weakly.  Again, by (\ref{e:13}) it suffices to
show that  $\sum_{j=2}^\infty (1-\sigma_j) =\infty$. The only
difference from the proof of part (i) is that the inequality used in
(\ref {e:15}) does no longer hold since $\delta_j >0$. However,
since $\delta_j\to0$,   we have, for $j$ large enough,
\[
1-\sigma_j=\frac{(1-\delta_{j-1}+\delta_j)(\delta_{j-1}-\delta_j)}{(1+\delta_j)(2\delta_{j-1}-\delta_j)}
> \frac{1}{2}\frac{\delta_{j-1}-\delta_j}{2\delta_{j-1}-\delta_j } >
\frac{1}{4}\frac{\delta_{j-1}-\delta_j}{\delta_{j-1}}>0.
\]
Then the same argument as in part (i) proves the claim.

\ep

The next special case is also based on iterated applications of Lemma \ref {L:2.1}  and shares part of the construction with the previous lemma, but with a different proof of the weak convergence of the vector sequence.

\bL{L:4.2}
Let $A= \sum_{i=1}^\infty (1+\mu_i) (e_i\otimes e_i)+ \sum_{i=1}^\infty (1-\lambda_i)(f_i \otimes f_i)$ where  $\{e_i, \, f_i\}$ are mutually orthogonal unit vectors, $\mu_i >0$, $0 < \lambda_i < 1$,  for all   $i$,  $\sum_{i=1}^\infty \lambda_i = \sum_{i=1}^\infty \mu_i< \infty$, and $\sum_{i=1}^m \lambda_i \neq \sum_{i=1}^n \mu_i$ for every $n,\, m\in \mathbb N$. Then $A$ is a strong sum of projections.
\eL
\bp

Since by hypothesis $\lambda_1\ne \mu_1$, we assume  that $\lambda_1>\mu_1$ and leave to the reader the  similar proof for the case when  $\lambda_1<\mu_1$.
 Since $\lambda_1 < \sum_{j=1}^\infty \mu_j $, there is a smallest integer $n_1$ for which $\lambda_1 < \sum_{j=1}^{n_1} \mu_j $. Similarly, there is  a smallest integer $m_1$  for which $ \sum_{j=1}^{m_1} \lambda_j > \sum_{j=1}^{n_1} \mu_j$. From here we obtain recursively the  strictly increasing integer  sequences  $ \{m_k\}$, $ \{n_k\}$, starting with $n_o=0,\, m_o=1$,  for which
\be{e:17}
\sum_{j=1}^{n_{k-1}} \mu_j
\le\sum_{j=1}^{n_k-1} \mu_j
< \sum_{j=1}^{m_{k-1}} \lambda_j
\le  \sum_{j=1}^{m_k-1} \lambda_j
< \sum_{j=1}^{n_k} \mu_j
\le \sum_{j=1}^{n_{k+1}-1} \mu_j
< \sum_{j=1}^{m_k} \lambda_j.
\ee
Set

\[
A_j:=
\begin{cases}
(1-\lambda_1)(f_1\otimes f_1) & j=1\\
(1-\lambda_{j-n_k})(f_{j-n_k}\otimes f_{j-n_k}) &m_{k-1}+n_k<j\le m_k+n_k\\
(1+\mu_{j-m_k})(e_{j-m_k}\otimes e_{j-m_k})&m_{k}+n_k<j\le m_k+n_{k+1}
\end{cases}
\]

Since $A$ is the sum of  two series  which  converge unconditionally, we can rearrange its summands to obtain  $A= \sum_{i=1}^\infty A_i$ (in the strong topology.)  Explicitly, for $j > 1$,
\[
\sum_{i=1}^j A_i=
\begin{cases}
\sum_{i=1}^{n_k}(1+\mu_i)(e_i \otimes e_i)+ \sum_{i=1}^{j-n_k}(1-\lambda_i)(f_i \otimes f_i)  &m_{k-1}+n_k\le j\le m_k+n_k\\
\sum_{i=1}^{j-m_k}(1+\mu_i)(e_i \otimes e_i) + \sum_{i=1}^{m_k}(1-\lambda_i)(f_i \otimes f_i)  &m_{k}+n_k\le j\le m_k+n_{k+1}.
\end{cases}
\]
Define
\[
\delta_j=
\begin{cases} -\lambda_1 \quad &j=1\\
\sum_{i=1}^{n_k}\mu_i - \sum_{i=1}^{j-n_k}\lambda_i &m_{k-1}+n_k< j\le m_k+n_k\\
\sum_{i=1}^{j-m_k}\mu_i - \sum_{i=1}^{m_k}\lambda_i &m_{k}+n_k <  j\le m_k+n_{k+1}.
\end{cases}
\]
Then from (\ref {e:17}) we have
\begin{alignat}{2}\label{e:18}
& \delta_j >0  \quad \quad &m_{k-1}+n_k\le j< m_k+n_k\\
& \delta_j <0  \quad \quad &m_{k}+n_k\le j< m_k+n_{k+1}\notag
\end{alignat}
and
\be{e:19}
\delta_j -\delta_{j-1}=
\begin{cases}
-\lambda_{j-n_k}< 0 &m_{k-1}+n_k<  j\le m_k+n_k\\
\phantom{a}\mu_{j-m_k} > 0&m_{k}+n_k<  j\le m_k+n_{k+1}.
\end{cases}
\ee
Thus
\begin{alignat}{2}\label{e:20}
\min\{\delta_j\mid m_{k-1}+n_k\le j\le m_k+n_{k+1}\} &= \delta_{n_k+m_k} \quad & (\text{by (\ref {e:19}}))\\
&= \delta_{n_k+m_k-1} - \lambda_{m_k} \quad & (\text{by (\ref {e:19}}))\notag\\
& > - \lambda_{m_k} & (\text{by (\ref {e:18}}))\notag\\
& > -1 &( \text{by  hypothesis})\notag.
\end{alignat}
Moreover,
\begin{alignat}{2}\label{e:21}
&2\delta_{j-1}-\delta_j = \delta_{j-1}+\lambda_{j-n_k} > \delta_{j-1}> 0\quad &m_{k-1}+n_k<  j\le m_k+n_k\\
&2\delta_{j-1}-\delta_j  = \delta_{j-1}- \mu_{j-m_k} <\delta_{j-1}<  0\quad &m_{k}+n_k<  j\le m_k+n_{k+1}.\notag
\end{alignat}

Define the sequence $\sigma_j$ as in (\ref {e:10}). From  (\ref {e:18}), (\ref {e:19}), and (\ref {e:21}), we see  that for every $j$, $\delta_{j-1}$,  $\delta_{j-1}-\delta_j $, and $2\delta_{j-1}-\delta_j$ have the same sign. Since furthermore $1+\delta_j> 0$ by (\ref {e:20}) and $1+\delta_j-\delta_{j-1}>0$ by \eqref {e:19} and \eqref {e:20},  it follows that $0 < \sigma_j< 1.$

Now let $J_k:=m_k+n_{k+1}$. Then we have by (\ref{e:18}) that $\delta_{J_k-1}<0<\delta_{J_k}$, by \eqref {e:20}  that $1+ \delta_{J_k}>0$, and hence   $\frac{1+\delta_{J_k-1} } { 1+\delta_{J_k}} < 1$. As  $2\delta_{J_k-1}-\delta_{J_k} <0$ by  (\ref {e:21}), we also have  $0 <   \frac { \delta_{J_k-1} } { 2\delta_{J_k-1}-\delta_j } < \2$
  and thus
 \be{e:22}
\sigma_{m_k+n_{k+1}} < \frac{1}{2}= \sigma_{J_k}.
\ee

Having concluded these preliminary computations, we define recursively the sequence of unit vectors
\be{e:23}
v_j =\begin{cases} f_1&j=1\\
 \sqrt{\sigma_j} v_{j-1} +  \sqrt{1-\sigma_j}  f_{j-n_k} &m_{k-1}+n_k<  j\le m_k+n_k\\
 \sqrt{\sigma_j} v_{j-1} +  \sqrt{1-\sigma_j} e_{j-m_k} &m_k+n_k<  j\le m_k+n_{k+1}.
\end{cases}
\ee
Notice that
\be{e:24}
v_{j}\in \begin{cases}
\text{span}\, \{f_1, \dots, f_{j-n_k}, e_1,\dots, e_{n_k}\}& m_{k-1}+n_k\le   j\le m_k+n_k\\
\text{span}\, \{f_1, \dots, f_{m_k}, e_1,\dots, e_{j-m_k}\}& m_{k}+n_k\le   j\le m_k+n_{k+1}.
\end{cases}
\ee
Now we claim that there is a sequence of rank-one projections $P_j$ for which
\be{e:25}
\sum_{j=1}^n A_j = \sum_{j=1}^{n-1}P_j + (1+\delta_n)(v_n\otimes v_n) \quad \text{for  }n \ge 2.
\ee
By Lemma \ref {L:2.1} 
\begin{alignat*}{2}
A_1&+A_2=(1-\lambda_1)(f_1 \otimes f_1)+(1+\mu_1)(e_1\otimes e_1)  & (\text{by definition})\\
&=  P_1 + (1+\mu_1-\lambda_1)(v\otimes v)\quad &(\text{by Lemma \ref {L:2.1} })\\
&= P_1 + (1+\delta_2)(v\otimes v) & (\text{by definition})
\end{alignat*}
where $P_1$ is a rank-one projection and  by  \eqref {e:3}, \eqref {e:2}, $v = \sqrt{\nu} f_1 + \sqrt{1-\nu } e_1$ and
\[
\nu =  \frac{(1-\lambda_1)\lambda_1}{(1+\mu_1-\lambda_1)(\mu_1+\lambda_1)} =\sigma_2
\]
and hence $v= v_2$.

Since $\delta_2 <0$ by \eqref {e:18} and  $v_2\perp e_2$, we can  apply Lemma \ref {L:2.1}  to 
\[
(1+\delta_2)(v_2\otimes v_2)  + (1+\mu_2)(e_2\otimes e_2)
\]
 and continue the process.  Assume  the construction up to  $j-1$, where  $m_{k-1}+n_k<  j\le m_k+n_k$ for some $k$.   Then

\ba
\sum_{i=1}^j A_i&= \sum_{i=1}^{j-1} A_i + (1-\lambda_{j-n_k})(f_{j-n_k}\otimes f_{j-n_k})\\
&= \sum_{1-1}^{j-2}P_i + (1+\delta_{j-1})(v_{j-1}\otimes v_{j-1})+ (1-\lambda_{j-n_k})(f_{j-n_k}\otimes f_{j-n_k}).
\end{align*}
Now $v_{j-1}\perp f_{j-n_k}$ by (\ref {e:24}) and $\delta_{j-1}>0$ by (\ref {e:18}), so we can apply Lemma \ref {L:2.1} and obtain
\[
(1+\delta_{j-1})(v_{j-1}\otimes v_{j-1})+ (1-\lambda_{j-n_k})(f_{j-n_k}\otimes f_{j-n_k})
= P_{j-1} + (1+\delta_{j-1}- \lambda_{j-n_k})(v\otimes v)
\]
where $P_{j-1}$ is a rank-one projection,  and  by  \eqref {e:3}, \eqref {e:2},  $v= \sqrt{\nu} f_{j-n_k} + \sqrt{1 - \nu} \, v_{j-1}$ and
\begin{alignat*}{2}
\nu &= \frac{(1- \lambda_{j-n_k})\lambda_{j-n_k}}{(1+\delta_{j-1}-\lambda_{j-n_k})(\delta_{j-1}+\lambda_{j-n_k})} &\quad \quad (\text {by Lemma \ref {L:2.1} } )\\
&=  \frac{(1+\delta_j-\delta_{j-1})(\delta_{j-1}- \delta_j)}{(1+\delta_{j})(2\delta_{j-1}- \delta_j)}& (\text{since $\delta_j = \delta_{j-1}-\lambda_{j-n_k}$ by (\ref {e:19})})\\
&=1-\sigma_j & (\text{by (\ref {e:11})})
\end{alignat*}
But then,  $v= v_j$ and since also $\delta_j = \delta_{j-1}-\lambda_{j-n_k}$, we see that
(\ref{e:25}) is satisfied for $j$.  We leave to the reader the similar proof for the case when  $m_{k}+n_k<  j\le m_k+n_{k+1}$ for some $k$.
We thus have for all $n$
\[
A-\sum_{i=1}^{n}P_i= (1+\delta_{n+1})(v_{n+1}\otimes v_{n+1})+ \sum_{i=n+1}^\infty A_i,
\]
Since  $ \sum_{i=n+1}^\infty A_i \underset{s}\to 0$, as in the proof of Lemma \ref {L:4.1},
in order to prove that $A = \sum_{i=1}^\infty P_i$ in the strong topology, it suffices to show that  the projections $v_j\otimes v_j\underset{s}\to 0$, or, equivalently, that the sequence of unit vectors $v_j\underset{w}\to 0$. Since all $v_j\in \text{span}\, \{f_i,e_i\}$, it suffices to prove that $(v_j, f_q)\to 0$ and $(v_j, e_q)\to 0$ for all $q$.

Fix $q\in \mathbb N$ and choose $h$ such that $m_h \ge q$ and $n_h \ge q$ and let $w=v_{m_h+n_h}$. From (\ref {e:23}) we have
\[
v_{m_h+n_h+1}- \sqrt{\sigma _{m_h+n_h+1} } w =  \sqrt{1-\sigma _{m_h+n_h+1}}e_{n_h+1}
 \in  \{f_1, \dots, f_{m_h}, e_1,\dots, e_{n_h} \}^\perp.
\]
Iterating,
\[
v_{j}- \Big(\prod_{i=m_h+n_h+1}^j\sqrt{\sigma _i }\Big)w \in \{f_1, \dots, f_{m_h}, e_1,\dots, e_{n_h} \}^\perp\quad\text{for every }j > m_h+n_h.
\]
In particular for every $j > m_h+n_h$,
\[
(v_j, f_q) = \Big(\prod _{i=m_h+n_h+1}^j \sqrt{\sigma_i }\Big)~(w, f_q)\quad\text{and}\quad
(v_j, e_q) = \Big(\prod _{i=m_h+n_h+1}^j \sqrt{\sigma_i }\Big)~(w, e_q).
\]
Since $0< \sigma_i < 1$  for all $i$ by (\ref {e:10}) and (\ref
{e:11}) and
 $\sigma_i < \2$ infinitely often by  (\ref
{e:22}), we see that  $\prod _{i=m_h+n_h+1}^j\sqrt{\sigma_i}\to0$
and hence $v_j\to 0$ weakly, which concludes the proof.

\ep
\bT{T:4.3} Let $A\in B(H)^+$ and assume that
$\tr(A_-)\le \tr(A_+)<\infty$  and  \linebreak $\tr(A_+)-\tr(A_-)  \in
\mathbb N\cup\{0\}$. Then $A$ is a strong sum of projections.
 \eT
\bp
Since $A_+$ and $A_-$ are of trace-class and supported in orthogonal subspaces,
they are simultaneously diagonalizable, so we can  set $A_-= \sum_{i=1}^M \lambda_i (f_i \otimes  f_i)$,
$A_+= \sum_{j=1}^N\mu_j (e_j \otimes e_j)$,  where $M, \, N \in \mathbb N \cup \{0\}\cup\{ \infty\}$,
 $\{f_i,\, e_j\}$ are mutually orthogonal unit vectors, and $0 < \lambda_i < 1$,
 $\mu_j >0$ or all $i$ and $j$. Let
\[
k:= \tr(A_+) -\tr(A_-)= \sum_{j=1}^N \mu_i-\sum_{i=1}^M \lambda_i.
\]
Since $\chi_A\{1\}$  is the sum of rank-one projections, we can by
(\ref {e:1}) assume without loss of generality that
\be{e:26} A=  \sum_{j=1}^N(1+\mu_j)(e_j \otimes e_j) +
\sum_{i=1}^M(1- \lambda_i) (f_i \otimes  f_i). \ee By the same proof
as in Lemma \ref {L:2.3}  we can decompose $A$ as the sum of $k$
rank-one projections and a positive operator $A'$ with $\tr(A'_+)
=\tr(A'_-)$. Thus we assume henceforth that $k=0$.

We need to consider four cases:

(a) when both $A_-$ and $A_+$ have finite rank (i.e.,  $N,M< \infty$),

(b) when $A_+$ has finite rank  and $A_-$ does not (i.e.,  $N<\infty$,
$M=\infty$),

(c) when $A_-$ has finite rank  and $A_+$ does not (i.e.,  $N=\infty$,
$M<\infty$), and

(d) when both have infinite rank (i.e.,  $N=M=\infty$.)

The case (a)   is given by Lemma  \ref {L:2.3}  (i).

Consider the case (b).  If $N>1$, choose an $m\in\mathbb N$ for
which 
\[
\sum_{j=1}^{N-1}\mu_j < \sum_{i=1}^m\lambda_i <
\sum_{j=1}^N\mu_j, \quad \text{i.e., } \quad 0< \sum_{j=1}^N \mu_j  - \sum_{i=1}^m\lambda_i <
\mu_N. 
\]
By Lemma \ref {L:2.3}  (ii)  there are $m+N$ rank-one
projections $P_k$ for which
\[
\sum_{j=1}^N(1+\mu_j)(e_j\otimes e_j)+\sum_{i=1}^m (1-\lambda_i)(f_i \otimes  f_i)=  \sum_{k=1}^{m+N-1}P_k + \Big(1+ \sum_{j=1}^N\mu_j -  \sum_{i=1}^m\lambda_i\Big)P_{m+N}
\]
Thus
\[
A':= A -  \sum_{k=1}^{m+N-1}P_k
= \Big (1+ \sum_{j=1}^N\mu_j -  \sum_{i=1}^m\lambda_i\Big) P_{m+N} + \sum_{j=m+1}^\infty(1-\lambda_j)(f_j\otimes f_j)
\]
Since $P_{m+N} \perp f_j$ for all $j>m$ and  $\sum_{j=1}^N\mu_j -  \sum_{i=1}^m\lambda_i =  \sum_{i=m+1}^\infty \lambda_i$,  we see that $A'$ satisfies the same conditions as $A$, but has $``N=1"$. Now  we obtain  by Lemma \ref {L:4.1} (ii) that $A'$ is a strong sum of projections and hence so is $A$.

The next case (c), when $N=\infty$ and $M<\infty$, is similar.  If
$M$ is not already $1$, choose an $n$ for which
$\sum_{j=1}^{n-1}\mu_j < \sum_{i=1}^{M-1}\lambda_i<
\sum_{j=1}^n\mu_j.$ Then, again by Lemma \ref {L:2.3}  (ii)
there are $M+n-1$ rank-one projections $P_k$ for which
\[
\sum_{j=1}^n(1+\mu_j)(e_j\otimes e_j)+ \sum_{i=1}^{M-1} (1-\lambda_i)(f_i \otimes  f_i)=
 \sum_{k=1}^{M+n-2}P_k + \Big(1+ \sum_{j=1}^n\mu_j -  \sum_{i=1}^{M-1}\lambda_i\Big)P_{M+n-1}
\]
Set
\[
A':= A -  \sum_{k=1}^{M+n-2}P_k
= (1-\lambda_M)(f_M\otimes f_M) + \Big(1+ \sum_{j=1}^n\mu_j -  \sum_{i=1}^{M-1}\lambda_i\Big)P_{M+n-1}
+ \sum_{j=n+1}^\infty(1+\mu_j)(e_j\otimes e_j).
\]
Since $P_{M+n-1}\le \sum_{j=1}^n(e_j\otimes e_j)+ \sum_{i=1}^{M-1} (f_i \otimes  f_i)$,
$P_{M+n-1}$ is orthogonal to the other rank-one summands of $A'$. Moreover,
\[
\lambda_M= \sum_{j=1}^n\mu_j -  \sum_{i=1}^{M-1}\lambda_i + \sum_{j=n+1}^\infty \mu_j.
\]
Now  we obtain  by Lemma \ref {L:4.1} (i) that $A'$ is a strong sum of projections and hence so is $A$.

In the last case (d),  both $N$ and $M$ are infinite  and
$\sum_{j=1}^\infty \lambda_j = \sum_{j=1}^\infty \mu_j$. Define
\[
\Phi_A:=\{(m,n)\in \mathbb N \times  \mathbb N \mid  \sum_{j=1}^m \lambda_j = \sum_{j=1}^n \mu_j\}.
\]
We need to treat the three possible cases separately, when $\Phi_A$ is infinite, when
it is finite and non-empty, and when it is empty.

If $\Phi_A$ is infinite, then rearrange it as $\Phi_A = \{(m_k, n_k)\}$ where the integer sequences
 $ \{m_k\}$, $ \{n_k\}$ are strictly increasing.  Set  $m_o=n_o = 0$ and by using the
 unconditional convergence of the series \eqref{e:26}, decompose $A$ as
\[
A=\sum_{k=1}^\infty \Big( \sum_{j=n_{k-1}+1}^{n_k}(1+ \mu _j) (e_j\otimes e_j)+
\sum_{j=m_{k-1}+1}^{m_k} (1-\lambda_j)(f_j\otimes f_j) \Big).
\]
 Since  $\sum_{j=n_{k-1}+1}^{n_k} \mu _j = \sum_{j=m_{k-1}+1}^{m_k} \lambda_j$, by
 Corollary \ref{C:2.5}   each summand 
 \[
 \sum_{j=n_{k-1}+1}^{n_k}(1+ \mu _j)
 (e_j\otimes e_j)+ \sum_{j=m_{k-1}+1}^{m_k} (1-\lambda_j)(f_j\otimes f_j) 
 \]
  is a sum of
 (finitely many) rank-one projections  and hence  $A$ is strong sum of projections.

If $\Phi_A$ is finite but not empty, it has a lexicographically largest element $(m,n)$ for which
$ \sum_{j=1}^{m }\lambda_j = \sum_{j=1}^{n }\mu_j$ but $ \sum_{j=1}^{m'} \lambda_j \ne \sum_{j=1}^{n'} \mu_j$
for any $m'>m$, $n'>n$. Now, again by Corollary \ref{C:2.5} ,
\[
 \sum_{j=1}^{n}(1+ \mu _j) (e_j\otimes e_j)+ \sum_{j=1}^m (1-\lambda_j)(f_j\otimes f_j )
\]
is the sum of rank-one projections and its remainder
\[
A':= A-  \sum_{j=1}^{n}(1+ \mu _j) (e_j\otimes e_j )+ \sum_{j=1}^{m} (1-\lambda_j)(f_j\otimes f_j)
\]
satisfies the same conditions as $A$, but in addition has
$\Phi_{A'}= \emptyset$. Finally, the crucial case when $\Phi_A=
\emptyset$ is given by Lemma \ref {L:4.2}.

\ep

In view of the necessary condition established in Theorem \ref {T:3.3}  (i), to conclude our study in $B(H)$  it remains to consider the case when $\tr(A_+)=\infty$. This will be done in Section \ref {S:6}.

\section{Type II factors: the  finite diagonalizable case}
In this section, we assume that $M$ is a type II factor with trace
$\tau$.  The following  key lemma is  also a consequence of Lemma
\ref {L:2.1} , or,  more precisely,  of Lemma \ref  {L:2.6}.

\bL{L:5.1} Let $A=
(1+\mu)E + (1-\lambda)F$  where $E$ and $F$ are finite projections,
$EF=0$, $\mu \ge 0$, $0\le \lambda \le 1$, and $\tau(A)\ge
\tau(R_A)$. Then $A$ is a strong sum of projections. \eL
 \bp 
 To avoid triviality, assume that $A\ne0$ and hence $E\ne 0$. The case when $\lambda = 1$
(resp. $\lambda=0$) is equivalent (resp., implied by) the case when
$F=0$, so assume that $0 < \lambda < 1$, and hence,  $R_A=E+F$. If
$\mu=0$,  then $ -\lambda\tau(F)=\tau(A)- \tau(R_A)\ge 0$, whence $
\lambda F = 0$, and then $A=E+F$ is already a projection. Thus
assume henceforth also that $\mu > 0$. Now consider first the key
case when $\tau(A)=\tau(R_A)$, i.e., $\mu\tau(E)= \lambda \tau(F)$.
In summary, assume that
\be {e:27} \mu > 0,\ \  0< \lambda < 1, \ \  \text{ and  } \ \
\mu\tau(E)= \lambda \tau(F)> 0. \ee 

If $\mu=\lambda$, then
$\tau(E)=\tau(F)$ and then $A$ is the sum of two  (equivalent)
projections by Lemma  \ref {L:2.6}.

 If $\mu<\lambda$, then $\tau(E)>\tau(F)$, and hence, there is some projection
 $E'\le E$ with \linebreak  $\tau(E')=\tau(F)$. Then $E'\sim F$ and by Lemma \ref {L:2.6} 
 there are projections $R_1,F_1 \in M$, $R_1,F_1 \le E'+F$ with
 $R_1\sim F_1\sim F$ for which 
 \[
 (1+\mu)E'+(1-\lambda)F= R_1 + (1+\mu-\lambda)F_1.
 \]
 Set $A_1:=A-R_1$, $E_1:=E-E'$, $\mu_1:=\mu$, and $\lambda_1:=\lambda- \mu$. Then
$E_1F_1=0$,  $E_1+ F_1\le E+F$, and
  \[
    \begin{cases}
    \mu_{1}=\mu >0\\
0 <  \lambda_{1}=\lambda-\mu < 1 \\
  \tau(E_{1})= \tau(E)-\tau(F) \\
   \tau(F_{1})= \tau(F)\\
    \end{cases}
  \]
Moreover, $A_1 = (1+\mu_1)E_1+ (1-\lambda_1)F_1$ and
 \[
 \mu_1  \tau(E_1) = \mu (\tau(E)- \tau(E'))= (\lambda-\mu)\tau(F)=\lambda_1\tau(F_1).
 \]
 Thus $A_1$ satisfies the same conditions (\ref{e:27}) as $A$ does.  
 
 Similarly, if  $\mu> \lambda$, and hence, $\tau(E)<\tau(F)$, choose a projection $F'\le F$ with \linebreak $\tau(F')=\tau(E)$. By the
same argument as above, there are projections  $R_1,E_1\in M,
\linebreak R_1, E_1\le E+F'$,  with
 $R_1\sim E_1\sim E$ for which 
 \[
 (1+\mu)E+(1-\lambda) F'= R_1 + (1+\mu-\lambda)E_1.
 \]

 Set $F_1:=F-F'$, $\mu_1:=\mu-\lambda$, $\lambda_1:=\lambda$, and $A_1:=A-R_1$. Then,
 again, $E_1F_1=0$,  $E_1+F_1\le E+F$, and
\[
 \begin{cases}\mu_{1}=\mu-\lambda > 0\\
 0 < \lambda_{1}=\lambda < 1 \\
  \tau(E_{1})= \tau(E)\\
   \tau(F_{1})= \tau(F)-\tau(E)
   \end{cases}
\]
Moreover, $A_1 = (1+\mu_1)E_1+ (1-\lambda_1)F_1$ and
$ \mu_1\tau(E_1) = \lambda_1\tau(F_1)$, i.e., here too $A_1$ satisfies the conditions (\ref{e:27}).

We can thus iterate the construction and find nonzero projections $E_k, F_k, R_k\in M$
with $E_kF_k=0$, $E_k +F_k \le E_{k-1}+F_{k-1}\le E+F$, positive operators $A_k= A_{k-1}-R_k$,
and scalars $\mu_k >0$ and $0< \lambda_k < 1$ for which
$\mu_k\tau(E_k) = \lambda_k\tau(F_k)$,  $A_k=(1+\mu_k)E_k+(1-\lambda_k)F_k$, and
\be{e:28}
 \begin{cases}\mu_{k+1}=\mu_k> 0\\
0<  \lambda_{k+1}=\lambda_k-\mu_k  <1\\
  \tau(E_{k+1})= \tau(E_k)-\tau(F_k) \\
   \tau(F_{k+1})= \tau(F_k)\\
    \end{cases} \quad \text{if } \mu_k < \lambda_k
    \ee
    \be{e:29}
   \begin{cases}\mu_{k+1}=\mu_k-\lambda_k> 0\\
0 < \lambda_{k+1}=\lambda_k < 1\\
  \tau(E_{k+1})= \tau(E_k)\\
   \tau(F_{k+1})= \tau(F_k)-\tau(E_k)x
   \end{cases}
   \quad  \text{if } \mu_k > \lambda_k
\ee Thus for every $k$, $A= \sum_{j=1}^kR_j + A_k$. This
construction terminates if for some $k$ we have
$\mu_{k}=\lambda_{k}$, in which case $A_{k}$ is the sum of two
projections, and hence, $A$ is the sum of $k+2$ projections. Thus
assume henceforth that $\mu_{k}\ne\lambda_{k}$ for every $k$.

By construction, both sequences $\tau(E_k)$ and $\tau(F_k)$ are
monotone non-increasing, and hence, both converge. Let $\alpha:=\lim \tau(E_k)
$ and $ \beta:= \lim \tau(F_k)$. The sequences $\mu_k$ and
$\lambda_k$ are also monotone non-increasing. If $\mu_{k+n}=\mu_k$
for some $k$ and $n\in \mathbb N$, then by (\ref {e:28}),
$\lambda_{k+n}=\lambda_k-n\mu_k$. Thus there must be a largest such
$n$, i.e., the sequence $\mu_k$ cannot be eventually constant and,
similarly, neither  can be the sequence $\lambda _k$.  Thus, both
inequalities $\mu_k< \lambda_k$ and $\mu_k> \lambda_k$ must occur
for infinitely many indices. Thus it follows from (\ref {e:28})
that $\alpha=\alpha-\beta$, and it follows from  (\ref {e:29})
that  $\beta = \alpha - \beta$, whence $\alpha=\beta=0$. As a
consequence, $||E_k||_1\to0$ and $||F_k||_1\to0$, and hence,
$E_k\underset{s}\to 0$ and $F_k\underset{s}\to 0$; this implication
is well known,  the reader is referred to \linebreak \cite[Exercise
8.7.39]{KrRj2}. Thus $A_k\underset{s}\to0$, and hence,
$A=\sum_{j=1}^\infty R_j $ where the convergence is also in the
strong operator topology.

We now consider the remaining case when $\mu>0$, $0< \lambda < 1$, and
$\mu\tau(E)>\lambda \tau(F)\ge 0$. Since $M$ is of type II, we can decompose
$E=E_1+E_2+E_3$ into the sum of three mutually orthogonal projections with
the following traces ($\lfloor{\mu}\rfloor$ denotes the integer part of $\mu$):
 \ba
 &\tau(E_1)= \frac{\lambda}{\mu}\tau(F)\ge 0 \\
 &\tau(E_2)= \frac{\mu\tau(E)-\lambda\tau(F)}{1+ \lfloor{\mu}\rfloor}> 0\\
  &\tau(E_3)=  \frac{(1-\mu+ \lfloor{\mu}\rfloor)(\mu\tau(E)-\lambda\tau(F))}{\mu(1+ \lfloor{\mu}\rfloor)}>0.
  \end{align*}
  Let
   \ba
&  A_1:= (1+\mu)E_1 + (1-\lambda)F \\
& A_2:= (\mu- \lfloor{\mu}\rfloor)E_2+ (1+\mu )E_3=  ((1-(1+ \lfloor{\mu}\rfloor-\mu)E_2+ (1+\mu )E_3\\
&A_3:=(1+\lfloor{\mu}\rfloor)E_2
     \end{align*}
Thus $A= A_1+A_2+A_3$. If $\tau(E_1)=0$, then $\lambda\tau(F)=0$,
and hence, $A_1= F$ is already a projection. If $\tau(E_1)\ne0$,
then $A_1$ satisfies the conditions of (\ref {e:27}), and hence, it
is a strong sum of projections. If $\mu\in \mathbb N$, then
$A_2=(1+\mu)E_3$ is the sum of $1+\mu$ projections. If $\mu\ne
\lfloor{\mu}\rfloor$, then it is easy to verify that also $A_2$
satisfies the conditions of (\ref {e:27}), and hence, is a strong sum
of projections. Finally,  $A_3$ is always trivially the sum of
$1+\lfloor{\mu}\rfloor $ projections, which concludes the proof. 

\ep

Now we consider positive  diagonalizable operators in $M$, namely,
those  operators of the form $A=\sum_k\gamma_kG_k$ where $G_k\in M$
are mutually orthogonal projections  and $\gamma_k> 0$, and the series, if infinite, converges in the strong operator topology.

\bT {T:5.2}
Let $M$ be a type II  factor with trace $\tau$ and let $A\in M$ be a positive
diagonalizable operator. If  $\tau(A_+)\ge \tau(A_-)$, then $A$ is a strong sum of projections.
\eT
\bp
To   avoid triviality,  assume that $A\ne 0$.  By renaming appropriately the coefficients and using the semifiniteness of
 $M$ to split the projections into sums of projections with finite trace, we rewrite  $A$ as
\[
A= \sum_{j=1}^N (1+ \mu_j)E_j+\sum_{i=1}^K (1-\lambda_i)F_i
\]
 with $N, K \in \mathbb N \cup \{0\}\cup \{\infty]$, $E_j, F_i$  mutually
 orthogonal finite projections, $\mu_j > 0$, and $0\le \lambda_i<1$ for all
 $j$ and $i$, and with the series converging strongly if $N$ or $K$  are infinite.
  Again, we use the convention that  if  $N$ or $K$ are zero then $A$ is the sum
of  only one series. Since  $\sum\{(1-\lambda_i)F_i \mid \lambda_i =0\}$
 is already a projection,
   we can further assume without loss of generality that $\lambda_i > 0$ for all $i$. Then
\[
A_+=  \sum_{j=1}^N \mu_j E_j ,\quad A_-= \sum_{i=1}^K \lambda_iF_i,
\quad \text{and hence,}\quad \sum_{j=1}^N \mu_j \tau(E_j)\ge
\sum_{j=1}^K \lambda_j \tau(F_j).
\]
In particular, $N > 0$. Assume $K>0$. Then $\lambda_1\tau(F_1)\le
\sum_{j=1}^N\mu_j\tau(E_j)$. Since $M$ is of type II, we can find
projections $E_{j1}\le E_j$ such that  $\lambda_1\tau(F_1)=
\sum_{j=1}^N\mu_j\tau(E_{j1})$. Then decompose $F_1=
\sum_{j=1}^NF_{j1}$ into mutually orthogonal projections so that
$\lambda_1\tau(F_{j1})=\mu_j\tau(E_{j1})$ for every $j$. If $K>1$ we
have
\[
 \sum_{i=2}^K \lambda_i \tau(F_i)\le \sum_{j=1}^N \mu_j
 \tau(E_j-E_{j1}),
\]
and hence, we can iterate the process. Thus  for every $i$ and $j$
we decompose $F_i = \sum_{j=1}^NF_{ji}$ into mutually orthogonal
projections and further  find mutually orthogonal projections \linebreak 
$E_{ji}\le E_j$ so that
 $\lambda_i\tau(F_{ji})=\mu_j\tau(E_{ji})$. Set $E_{jo}:= E_j-\sum_{i=1}^K E_{ji}$.
 Then
 \[
 A= \sum_{j=1}^N \sum_{i=1}^K  \Big( (1+\mu_j)E_{ji}+(1- \lambda_i)F_{ji}   \Big) + \sum_{j=1}^N (1+\mu_j) E_{jo}.
 \]
By Lemma \ref {L:5.1}, each summand $(1+\mu_j)E_{ji}+(1-
\lambda_i)F_{ji}$ and  $(1+ \mu_j)E_{jo}$ is a strong sum of
projections, and hence, so is $A$. In the case that $K=0$,
$A=\sum_{j=1}^N (1+\mu_j) E_{j}, $ and hence, it is also the strong
sum of projections by the same reasoning.

\ep
As the following examples show, the condition that $A$ is diagonalizable is not necessary for $A$
to be a strong sum of projections.

\begin{example}\label {E:5.3}
Let $M$ be a type $II_1$ factor, let   $P\in M$ be a projections with $P \sim P^\perp$,
let $\mathscr A$ (resp.,  $\mathscr B$) be a masa in $M_P$  (resp., in $M_{P^\perp}$).
By properly scaling the spectral resolution of  a generator of $\mathscr A$ we can find 
a monotone increasing strongly continuous net of projections $\{E_t\}_{t\in[0,\2]}$ in
$\mathscr A$ with $\tau(E_t)=t$.
\item [(i)] Assume that $\mathscr A$ and $\mathscr B$ are conjugate in $M$, and hence there is a 
selfadjoint unitary  $U\in M$  for which $U\mathscr AU^*=\mathscr B$. Define
\[
E_t:=I - UE_{1-t}U^* ~~\text{for }t\in (\2,1]\quad\text{and}\quad A= \int_0^{\2}(1+t)dE_t+\int_{\2}^1tdE_t.
\]
Then $\{E_t\}_{t\in[0,1]}$ is flag, namely  a monotone increasing strongly continuous net of projections with $\tau(E_t)=t$ for all $t\in [0,1]$ and $A$ is not diagonalizable (in fact, it has no eigenvalues). Furthermore,  $0\le A\le \frac{3}{2}I\le 2I$, $R_A=I$, 
and it is easy to verify that
\ba
UAU^*&=  \int_0^{\2}(1+t)d(UE_tU^*)+ \int_{\2}^1td(UE_tU^*)\\
&=\int_0^{\2}(1+t)d(I- E_{1-t})+ \int_{\2}^1td(I-E_{1-t})\\
&= \int_{\2}^1(2-t)dE_t+\int_0^{\2}(2-(1+t))dE_t\\
&=2I-A.
\end{align*}
Thus by Proposition \ref {L:2.1} , $A$ is the sum of two projections.
\item[(ii)] Assume that  $\mathscr A$ and $\mathscr B$ are not conjugate in $M$. Such a case can be easily obtained by choosing $P$ so that $M\sim M_P\sim M_{P^\perp}$, choosing two non-conjugate masas $\mathscr A_o$ and $\mathscr B_o$ in $M$ (e.g., a Cartan masa and a singular one) and defining $\mathscr A$ and $\mathscr B$ to be the compressions of $\mathscr A_o$ and $\mathscr B_o$ to $M_P$ and $M_{P^\perp}$ respectively. Complete  $\{E_t\}_{t\in[0,\2]}$ to be a flag in $M$ by defining $E_t:=P+ F_t ~~\text{for }t\in (\2,1]$ where $\{F_t\}_{t\in[\2,1]}$ is an arbitrary monotone increasing strongly continuous net of projections in $\mathscr B$ with $\tau(F_t)= t-\2$. Define as in (i)
\[
A= \int_0^{\2}(1+t)dE_t+\int_{\2}^1tdE_t.
\]
Again, $A$ is not diagonalizable, in fact it has no,  $0\le A\le \frac{3}{2}I\le 2I$, and $R_A=I$. By  Proposition \ref {P:2.10}  , $A$ is the sum of two projections if and only if $UAU^*=2I-A$ for some unitary $U\in M$. Reasoning as in the proof of Proposition \ref {P:2.10} , it is simple to see that if such a unitary existed, we would have $\mathscr B= U\mathscr AU^*$ and hence $\mathscr A$ and $\mathscr B$ would be conjugate, against the assumption. Thus $A$ cannot be the sum of two projections in $M$. However we do not know whether $A$ is a strong sum of projections in $M$ or not.
\end {example}

\begin{question}\label{Q:5.4} Can the condition that $A$ is diagonalizable be removed from Theorem \ref {T:5.2}?
\end{question}

\section{The infinite case}\label{S:6}
In this section we assume that $M$ is an infinite factor, i.e., of  type I$_\infty$,  type II$_\infty$, or type III. The following lemma is the key to the proof of Theorem \ref {T:1.1} in this case.

\bL{L:6.1} Let $A=\sum_{j=1}^\infty(1+\mu_j) E_j+ (1-\lambda) F $  where $\{E_j, F\}$ are  mutually orthogonal equivalent projections in $M,~\mu_j >0$, $0< \lambda_j\le 1$, and $\sup \mu_j < \infty$. If $\sum _{j=1}^\infty\mu_j=\infty$, then $A$ is a strong sum of projections in $M$.
\eL

\bp
Let $n_1\ge 1$ be the smallest integer for which $\sum_{j=1}^{n_1}\mu_j \ge \lambda $. Such an integer exists because $\sum_{j=1}^\infty \mu_j=\infty$. Set
\[
\alpha_1:= \begin{cases} \lambda&\text{if } n_1=1\\
\lambda -  \sum_{j=1}^{n_1-1}\mu_j&\text{if } n_1>1
\end{cases}\quad \text{and}\quad 1-\beta_1:=  \mu_{n_1}-\alpha_1 - \lfloor\mu_{n_1}-\alpha_1\rfloor
\]
where  $\lfloor x \rfloor$ denotes the integer part of $x$.
Then $\mu_{n_1}\ge \alpha_1$, $0<  \alpha_1\le 1$, and  $0< \beta_1 \le 1 $. The positive operator
\[
D_1:=\sum_{j=1}^{n_1-1}(1+\mu_j)E_j + (1+\alpha_1 +  \lfloor\mu_{n_1}-\alpha_1\rfloor)E_{n_1} + (1-\lambda) F
\]
is a linear combination of  $n:=\begin{cases}n_1+1 &\text{if } \lambda\ne1\\
n_1&\text{if } \lambda =1
\end{cases} $ mutually orthogonal equivalent projections in $M$ and the sum of their coefficients is
 $k_1:=n_1+1+  \lfloor\mu_{n_1}-\alpha_1\rfloor.$  Since $k_1\in \mathbb N$ and $k_1\ge n$, by Lemma \ref {L:2.6}  (ii), $D_1$ is the sum of $k_1$ (equivalent) projections.

Next, we apply the same construction to  the  ``remainder"
\[
A-D_1=   \sum_{j=n_1+1}^\infty(1+\mu_j)E_j   + (1-\beta_1)E_{n_1}
\]
where now $\beta_1$ plays the role of $\lambda$ and $E_{n_1}$ the role of $F$. Iterating we find
an increasing sequence of indices $n_k$ and two sequences of positive numbers $0< \alpha_k,  \beta_k\le 1$
with $\mu_{n_k}\ge \alpha_k$ and  $1-\beta_k=\mu_{n_k}- \alpha_k - \lfloor
\mu_{n_k}- \alpha_k\rfloor$. Then the positive operator
\[
D_k:= \sum_{j=n_{k-1}+1}^{n_k-1}(1+\mu_j)E_j +  (1+\alpha_k +  \lfloor\mu_{n_k}-\alpha_k\rfloor)
E_{n_k}+ (1-\beta_{k-1})E_{n_{k-1}}
 \]
is by Lemma \ref {L:2.6}  the sum of finitely many (equivalent) projections. But then
\be {e:30}
A-\sum_{j=1}^kD_j = \sum_{j=n_{k}+1}^\infty(1+\mu_j)E_j +(1-\beta_{k})E_{n_k}   \underset{s}{\to}0
\ee
because the projections $\{E_j\}$ are mutually orthogonal. Thus $A= \sum_{j=1}^{\infty}D_j $, where the series converges in the strong operator topology. Since each $D_k$ is the sum of  projections, so is $A$.

\ep

If $M$ is of type I and all projections $E_j$ and $F$ have rank-one, then we can  relax the
condition that they are mutually orthogonal. Indeed, orthogonality is not necessary to conclude
that each positive finite rank operator $D_k$ is the sum of  projections  (see Corollary
 \ref {L:2.6}  and also Lemma \ref {L:2.2} ), and assuming strong convergence of
 the series $\sum_{j=1}^\infty(1+\mu_j) E_j$ is sufficient to guarantee that
 $ \sum_{j=n_{k}+1}^\infty(1+\mu_j)E_j +(1-\beta_{k})E_{n_k}   \underset{s}{\to}0$ in
 \eqref{e:30}. Thus we have:

\bL{L:6.2}Let $A=\sum_{j=1}^\infty(1+\mu_j) E_j+ (1-\lambda) F $  where $E_j,\,F\in B(H)$ are
 rank-one projections, $\mu_j >0$, $0< \lambda_j\le 1$, and $\sum_{j=1}^\infty(1+\mu_j) E_j$ converges
 in the strong operator topology. If $\sum _{j=1}^\infty\mu_j=\infty$, then $A$ is a strong sum of projections.
\eL

\bP{P:6.3}Let   $A\in M^+$ and assume that there is some $\mu>0$ for which the spectral projection $\chi_A[1+\mu,\infty)$
is infinite. Then $A$ is  a strong sum of projections.
\eP
\bp
Let $E:=\chi_A[1+\mu,\infty)$, $B:=A-(1+\mu)E$, and let $\mathscr A$ be a masa  containing $A$. Then
 $B\in \mathscr A$. By \cite[Corollary 2.23]{SsZl}, $B$ can be decomposed into a norm
 converging series $B=\sum_{i=1}^\infty (1-\lambda_i)Q_i$ with $0\le \lambda_i<1$ and with
 the projections $Q_i\in\mathscr A$. (In fact  we can choose $1-\lambda_i = 2^{-i}$, but we do not need this fact here.) Some or all of  the projections $Q_i$ can be zero.
Since $M$ is infinite and $E\in \mathscr A$, by
 \cite[Theorem 3.18] {Kr84} (see also \cite [Corollary~31]  {Kv91}), we can decompose
 $E=\sum_{i=1}^\infty E_i$  into a sum of infinite projections $E_i\in\mathscr A$. Let \linebreak$A_i:= (1+\mu)E_i+(1-\lambda_i)Q_i$. Then $A=\sum_{i=1}^\infty A_i$. Thus it
  suffices to prove that $A_i$ is a strong sum of projections  for each $i$. Using
  the fact that $E_i,~Q_i\in\mathscr A$, and hence, they commute, it follows that $A_i$ is
  diagonalizable as
\[
A_i= (1+\mu)(E_i-E_iQ_i)+  (2+\mu-\lambda_i)E_iQ_i+(1-\lambda_i)(Q_i-E_iQ_i).
\]
Since $E_i$ is infinite, at least one of the two orthogonal projections $E_i-E_iQ_i$
 and $E_iQ_i$ must be infinite. Assume that $E_i-E_iQ_i$  is infinite. If  $Q_i= 0$, then $A_i= (1+\mu)E_i$ and the conclusion follows from Lemma \ref {L:6.1} by further decomposing  $E_i$ into a sum of infinitely many mutually orthogonal equivalent projections. 
 
 Thus assume that $Q_i\ne 0$ and  decompose  $2+\mu-\lambda_i= \sum_{n=1}^m(1-\gamma_n)$ into the sum of finitely many numbers
 $0<1- \gamma_n< 1$.  Next, 
decompose $E_i-E_iQ_i= \sum_{n=1}^{m+1}E^{(n)}_i$  into the sum of $m+1$  mutually orthogonal equivalent (infinite) projections $E^{(n)}_i$.  Then  further decompose each projection $E^{(n)}_i$ into a sum
of infinitely many mutually orthogonal projections $E_{ij}^{(n)}$ with
\[
E_{ij}^{(n)}\sim \begin{cases} E_iQ_i\quad&\text{for } 1\le n \le m\\
Q_i-E_iQ_i &\text{for } n=m+1.
\end{cases}
\]
Thus $E_i-E_iQ_i=\sum _{n=1}^{m+1}\sum_{j=1}^\infty E_{ij}^{(n)}$.
Define
\[
B_{i}^{(n)}= \begin{cases} \sum_{j=1}^\infty (1+\mu)E_{ij}^{(n)} + (1-\gamma_n)E_iQ_i
\quad&\text{for } 1\le n \le m\\
\sum_{j=1}^\infty (1+\mu)E_{ij}^{(m+1)}
+(1-\lambda_i)(Q_i-E_iQ_i)&\text{for } n=m+1.
\end{cases}
\]

By  construction,  $A_i= \sum _{n=1}^{m+1}B_{i}^{(n)}$. By Lemma \ref {L:6.1} ,
 all the operators $B_{i}^{(n)}$ are strong sums of projections and hence so is $A_i$.
 Finally, the case when $E_iQ_i$ is infinite is  similar and is left to the reader.

 \ep

An immediate consequence of this proposition is the sufficient condition in Theorem \ref {T:1.1} (iii) for the type III case.

\bC{C:6.4} Let $M$ be a type III factor,  $A\in M^+$, and 
either $A$ be a projection or $A$ satisfy $||A||>1$. Then $A$ is a strong
sum of projections.

 \eC \bp If  $||A||>1$, then there is some
$\mu>0$ for which the spectral projection $\chi_A[1+\mu,\infty)$
is nonzero and hence infinite. Then $A$ is a strong sum of
projections by Proposition \ref {P:6.3}. 

\ep

\bR{R:6.5}
\item [(i)] The condition that $\chi_A[1+\mu,\infty)$ is infinite for some $\mu>0$ is equivalent
 to the condition $||A||_{ess}>1$ where $||A||_{ess}$ is the \textit{essential norm}, i.e., the norm
  in the quotient $M/K$, where $K$ is the norm closed ideal generated by the finite projections of $M$.

  If $M=B(H)$, then $K$ is the ideal of compact operators $K(H)$ on $H$ and Proposition \ref {P:6.3}  provides another proof of \cite [Theorem 2 ]     {DFKLOW} stating that if $||A||_{ess}>1$, then $A$ is a strong sum of projections.
  
   If  $M$ is of type II$_\infty$,
   $K$ is the ideal of compact operators relative to $M$ introduced by Sonis \cite{Sm} and  Breuer \cite{Bm}
   (see also \cite{Kv77}). 
   
   If $M$ is of type III,  then $K=\{0\}$ and  $||A||_{ess}=  ||A||$.  
   \item[(ii)] 

If $M$ is semifinite and $A\in K^+$ is  a strong sum of projections then $\tau(R_A) < \infty$.
\eR
\bp
\item[(ii)]It is well known that $\tau(\chi_A(\gamma, \infty))< \infty$ for every $\gamma >0$ (e.g., see
\cite [Theorem 1.3] {Kv77}). In particular, $\tau(\chi_A(1, \infty))< \infty$, whence
 $\tau(A_+)< \infty$. Thus it follows from Theorem \ref {T:3.3}  that  $\tau(A_-)< \infty$.
  But $A_-\ge \2\chi_A(0, \2]$, hence $\tau(\chi_A(0, \2))< \infty$. From this follows that  $\tau(\chi_A(0, \infty))
   = \tau(\chi_A(0, \2])+ \tau(\chi_A( \2, \infty))< \infty$
\ep
 An alternative proof of (ii) for the case when $M=B(H)$ is that a strongly converging series of rank-one projections that converges to a compact operator must converge uniformly and hence be finite.

We can now prove the last part of the sufficiency in Theorem \ref {T:1.1}.

\bT{T:6.6} Let $M$ be type I$_\infty$ or type II$_\infty$. If
$\tau(A_+)= \infty$. Then $A$ is a strong sum of projections. \eT

\bp
By Proposition \ref {P:6.3}, we  just need to consider the case when
 $\chi_A[1+\mu, \infty)$ is finite for every $\mu>0$. Let $\mathscr A$ be a masa
 of $M$ containing $A$. Let $E'_1:=\chi_A[2, \infty)= \chi_A[2, ||A||]$ and
 $E'_j:=\chi_A[1+\frac{1}{j}, 1+\frac{1}{j-1})$ for  $j >1$. Then $\tau(E'_j)< \infty $ for all $j$. Since
\[
\sum_{j=1}^\infty \frac{1}{j}E'_j\le A_+\le (||A||-1)E'_1+ \sum_{j=2}^\infty \frac{1}{j-1}E'_j,
\]
we see that
$
(||A||-1)\tau(E'_1)+ \sum_{j=2}^\infty \frac{1}{j-1}\tau(E'_j)= \infty .
$
Then  also $\sum_{j=1}^\infty \frac{1}{j}\tau(E'_j) = \infty.$
Furthermore, $A\chi_A(1, ||A||] - \sum_{j=1}^\infty
(1+\frac{1}{j})E'_j\in \mathscr A^+$, $A\chi_A[0,1] \in \mathscr
A^+$, and 
\[
 A= \sum_{j=1}^\infty (1+\frac{1}{j})E'_j +
A\chi_A[0,1] + \Big(A\chi_A(1, ||A||] - \sum_{j=1}^\infty
(1+\frac{1}{j})E'_j\Big). 
\]
 Now we consider separately the case
when $M=B(H)$ and when $M$ is of type II. 

If $M=B(H)$, first
decompose by  \cite [Corollary 2.23]{SsZl} the positive operator

\[
B:= A\chi_A[0,1] + \Big(A\chi_A(1, ||A||] - \sum_{j=1}^\infty
(1+\frac{1}{j})E'_j\Big)
\]
 into a norm converging series $B= \sum_{i=1}^\infty (1-\lambda'_i)Q'_i$ with $0< \lambda'_i\le 1$ and with
 the projections $Q'_i\in\mathscr A$.  Some or all of  the projections $Q'_i$ can be zero.
Then,  further decompose the projections $E'_j$ and $Q'_i$ into
rank-one projections. Relabel the ensuing sequence of coefficients
$1+ \frac{1}{j}$ (resp., $1-\lambda'_i$) repeated
according to the multiplicity of the projections as $1+\mu_j$ (resp.
$1-\lambda_i)$. To take into account the case when there are only
finitely many non-zero projections $Q'_i$ and they all have finite
rank, allow $ \lambda_i=1$. Thus
\[
A= \sum_{j=1}^\infty (1+\mu_j)E_j+ \sum_{i=1}^\infty(1-\lambda_i)Q_i
\]
 where all the projections $E_j$ and $Q_i$ have rank-one, $\mu_j > 0, ~0 < \lambda_i\le 1$ for all $i$ and $j$, both series converge in the strong operator topology,  and $\sum_{j=1}^\infty \mu_j= \infty$.  Now further decompose $\mathbb N = \cup _{i=1}^\infty \Lambda_i$ into infinite disjoint subsets $\Lambda_i$ so that  for each $i$,  $\sum_{j\in \Lambda_i}\mu_j= \infty$. Then
 \[
 A= \sum_{ i=1 }^\infty \Bigg( \sum_{ j\in \Lambda_i } ( 1+\mu_j ) E_j+ (1-\lambda_i) Q_i\Bigg)
 \]
and each summand being the strong sum of projection by Lemma \ref {L:6.2}, so is $A$.

 Now assume that $M$ is of type II and again using  \cite [Corollary 2.23]{SsZl} decompose separately
$A\chi_A[0,1]$ and $A\chi_A(1, ||A||] - \sum_{j=1}^\infty (1+\frac{1}{j})E_j$ into two norm converging series of scalar multiples $1-\lambda'_i$ of projections $Q'_i\in \mathscr A$. If a projection $Q'_i$ in the series decomposing $A\chi_A[0,1]$ is infinite, by the semifiniteness of $M$ we can further decompose it into a strongly converging sum of mutually orthogonal finite projections in $M$. These projections are not necessarily in $\mathscr A$,  however,  being majorized by $\chi_A[0,1]$, they  are all orthogonal to and hence commute with all the projections $E'_j$.

Every projection $Q'_i$ in the series decomposing $A\chi_A(1, ||A||] - \sum_{j=1}^\infty (1+\frac{1}{j})E_j$ is in $\mathscr A$, is majorized by $\chi_A(1, ||A||]= \sum_{j=1}^\infty E_j$, and hence is the sum $Q'_i=\sum_{j=1}^\infty Q'_iE_j$ of finite projections $Q'_iE'_j$ which belong to $\mathscr A$ and hence commute with all the projections $E'_j$. Therefore,
\[
A= \sum_{j=1}^\infty (1+\frac{1}{j})E'_j+ \sum_{i=1}^\infty(1-\lambda_i)Q_i
\]
where for all $i$, $0\le \lambda_i< 1$ and $Q_i$ are finite projections that commute with each $E_j$.  Since $\sum_{j=1}^\infty \frac{1}{j}\tau(E'_j) = \infty$, we can choose an increasing sequence of indices $n_i$ for which  $\sum_{j=n_i+1}^{n_{i+1}}\frac{1}{j}\tau(E'_j)\ge \lambda_i\tau(Q_i)$ and let  $A_i:=   \sum_{j=n_i+1}^{n_{i+1}}(1+\frac{1}{j})E'_j+(1- \lambda_i)Q_i$. Then  $A= \sum_{i=1}^\infty A_i$ and

 \[
 \tau(A_i)=  \sum_{j=n_i+1}^{n_{i+1}}\tau(E'_j)+\tau(Q_i) + \sum_{j=n_i+1}^{n_{i+1}}\frac{1}{j}\tau(E'_j)- \lambda_i\tau(Q_i) \ge \tau(\big( \sum_{j=n_i+1}^{n_{i+1}}E'_j \big) \vee Q_i\big)= \tau(R_{A_i})
  \]
Since  $A_i$ is diagonalizable because $Q_i$ commutes with all $E'_j$, $A_i$ is a strong sum of projections by Theorem \ref {T:5.2} and hence so is $A$.
\ep

\end{document}